%% file: main.tex
\documentclass[12pt]{article}
\usepackage{geometry}
\usepackage[utf8]{inputenc}
\usepackage[T1]{fontenc}

\usepackage{amsmath,amsthm,amssymb}
\usepackage{graphicx}
\usepackage{epsfig}
\usepackage{color}

\usepackage[labelfont=bf]{caption}
\usepackage{subfigure}

\DeclareMathOperator{\atantwo}{atan2}

\def \R {{\rm I \mskip -3.0 mu R}} 
\def \E {{\rm I \mskip -3.0 mu E}}
\def\a{{\bf a}}
\def\b{{\bf b}}
\def\c{{\bf c}}
\def\C{{\bf C}}
\def\d{{\bf d}}

\def\i{{\bf i}}

\def\o{{\bf o}}
\def\p{{\bf p}}
\def\q{{\bf q}}
\def\r{{\bf r}}

\def\u{{\bf u}}
\def\vv{{\bf v}}
\def\x{{\bf x}}

\def\w{{\bf w}}
\def\X{{\bf X}}
\def\A{{\cal A}}
\def\B{{\cal B}}
\def\G{{\cal G}}
\def\U{{\cal U}}

\def\AMM{{\it Amer.\ Math.\ Monthly\ }}
\def\ACM{{\it Adv.\ Comp.\ Math.\ }}
\def\ACMTMS{{\it ACM Trans.\ Math.\ Software\ }}
\def\ACMTOG{{\it ACM Trans.\ Graphics\ }}
\def\CAD{{\it Comput.\ Aided Design }}
\def\CAEJ{{\it Comput.\ Aided Eng. J.\ }}
\def\CAGD{{\it Comput.\ Aided Geom.\ Design }}
\def\CG{{\it Computers \& Graphics }}
\def\CVGIP{{\it Comput.\ Vision, Graphics, Image\ Proc.\ }}
\def\IBMJRD{{\it IBM J.\ Res.\ Develop.\ }}
\def\JMAA{{\it J.\ Math.\ Anal.\ Appl.\ }}
\def\MC{{\it Math.\ Comp.\ }}
\def\MMAS{{\it Math.\ Methods\ Appl.\ Sci.\ }}
\def\NA{{\it Numer.\ Algor.\ }}
\def\PAMS{{\it Proc.\ Amer.\ Math.\ Soc.\ }}
\def\SIAMJNA{{\it SIAM J.\ Numer.\ Anal.\ }}
\def\SIAMR{{\it SIAM Rev.\ }}

\newcommand{\be}{\begin{equation}}
\newcommand{\ee}{\end{equation}}
\newcommand{\ba}{\begin{eqnarray}}
\newcommand{\ea}{\end{eqnarray}}
\newcommand{\bi}{\begin{itemize}}
\newcommand{\ei}{\end{itemize}}

\newtheorem{prn}{Proposition}

\newcommand{\prf}{\noindent{{\bf Proof} :\ }}
\newcommand{\QED}{\vrule height 1.4ex width 1.0ex depth -.1ex\ \medskip}

\begin{document}

\title{A local $C^2$ Hermite interpolation scheme with PH quintic splines for 3D data streams}

\author{
Carlotta Giannelli, Lorenzo Sacco, Alessandra Sestini\\
Dipartimento di Matematica e Informatica ``U. Dini,'' \\
Universit\`a degli Studi di Firenze, Viale Morgagni 67/A, I--50134 Firenze, Italy
}

\date{}

\maketitle

\thispagestyle{empty}

\begin{abstract}
The construction of  smooth spatial paths with Pythagorean-hodograph (PH) quintic spline biarcs is proposed. To facilitate real-time computations of $C^2$ PH quintic splines, an efficient \emph{local} data stream interpolation algorithm is introduced. Each  spline segment interpolates second and first order Hermite data at the 
initial and final end-point, respectively. In the spline extension of the scheme a $C^2$ smooth connection between successive spline segments is obtained by taking the locally required second-order derivative information from the previous segment.
Consequently, the data stream spline interpolant is globally $C^2$ continuous and can be constructed for arbitrary $C^1$ Hermite data configurations. A simple and effective selection of the free parameters that arise in the interpolation problem is proposed. The developed theoretical analysis proves the fourth approximation order of the local scheme while a selection of numerical examples confirms the same accuracy of its spline extension. In addition, the performances of the algorithm are also validated by considering its application to point stream interpolation with automatically generated first-order derivative information.
\end{abstract}

\thispagestyle{empty}

 \centerline{{\bf Keywords}:  Pythagorean-hodograph 
curves; Biarcs; Data stream interpolation;} 
\centerline{Hermite interpolation; Quaternions}

\bigskip
\centerline{e--mail: carlotta.giannelli@unifi.it,}
\centerline{lorenzo.sacco@unifi.it, alessandra.sestini@unifi.it}

\newpage

\setcounter{page}{1}

\thispagestyle{plain}

\input{section01}

%
\input{section02}

\input{section03}

\input{section04}

\input{section05}

\section{Closure}
\label{sec:closure}
Quintic Pythagorean-hodograph biarcs are here proposed to construct spatial $C^2$ spline interpolants. The locality of the scheme is suitable for real time interpolation of Hermite data streams.
A robust and effective data-dependent strategy to fix the three free parameters associated to each spline segment is proposed.
This strategy endows the scheme with fourth approximation order, have an easy implementation, and produce paths with suitable fair shape when non asymptotic data are considered. 
An application to point stream (without derivative information) interpolation is also presented by properly combining the proposed scheme with simple local formulas for derivative approximation.

\input{Appendix}


\def\AAA{{\it Abstr.\ Appl.\ Anal.\ }}
\def\AAAS{{\it Acta\ Aero.\ Astro.\ Sinica\ }}
\def\AACA{{\it Adv.\ Appl.\ Clifford\ Alg.\ }}
\def\AC{{\it Acta\ Cybernetica\ }}
\def\ACM{{\it Adv.\ Comp.\ Math.\ }}
\def\ACMTMS{{\it ACM Trans.\ Math.\ Software\ }}
\def\ACMTOG{{\it ACM Trans.\ Graphics\ }}
\def\AMC{{\it Appl.\ Math.\ Comp.\ }}
\def\AMM{{\it Amer.\ Math.\ Monthly\ }}
\def\ASMEJDSMC{{\it ASME J.\ Dyn.\ Syst.\ Meas.\ Control\ }}
\def\CAD{{\it Comput.\ Aided Design }}
\def\CAEJ{{\it Comput.\ Aided Eng. J.\ }}
\def\CAGD{{\it Comput.\ Aided Geom.\ Design }}
\def\CAM{{\it Comput.\ Appl.\ Math.\ }}
\def\CAVW{{\it Comput.\ Anim.\ Virt.\ Worlds }}
\def\CG{{\it Computers \& Graphics\ }}
\def\CMC{{\it Comput.\ Meas.\ Control\ }}
\def\CVGIP{{\it Comput.\ Vision, Graphics, Image\ Proc.\ }}
\def\EO{{\it Eng.\ Optim.\ }}
\def\GM{{\it Graph.\ Models\ }}
\def\IBMJRD{{\it IBM J.\ Res.\ Develop.\ }}
\def\IJAMT{{\it Int.\ J.\ Adv.\ Manufac.\ Tech.\ }}
\def\IJAT{{\it Int.\ J.\ Autom.\ Tech.\ }}
\def\IJCM{{\it Int.\ J.\ Comput.\ Math.\ }}
\def\IJMTM{{\it Int.\ J.\ Mach.\ Tools Manufac.\ }}
\def\IJPEM{{\it Int.\ J.\ Prec.\ Eng.\ Manufac.\ }}
\def\IJSSIS{{\it Int.\ J.\ Smart\ Sens.\ Intell.\ Syst.\ }}
\def\IMAJNA{{\it IMA J.\ Numer.\ Anal.\ }}
\def\JAE{{\it J.\ Aerosp.\ Eng.\ }} 
\def\JBCS{{\it J.\ Braz.\ Comput.\ Soc.\ }}
\def\JCADCG{{\it J.\ Comput.\ Aided Design \& Comput.\ Graphics }}
\def\JCAM{{\it J.\ Comput.\ Appl.\ Math.\ }}
\def\JGCD{{\it J.\ Guid.\ Contr.\ Dyn.\ }}
\def\JGG{{\it J.\ Geom.\ Graphics }}
\def\JIRS{{\it J.\ Intell.\ Robot\ Syst.\ }}
\def\JMAA{{\it J.\ Math.\ Anal.\ Appl.\ }}
\def\JNM{{\it J.\ Numer.\ Math.\ }}
\def\JSC{{\it J.\ Symb.\ Comput.\ }}
\def\MC{{\it Math.\ Comp.\ }}
\def\MCS{{\it Math.\ Comput.\ Simul.\ }}
\def\MIC{{\it Model.\ Ident.\ Control\ }}
\def\MMAS{{\it Math.\ Methods\ Appl.\ Sci.\ }}
\def\NA{{\it Numer.\ Algor.\ }}
\def\NMTMA{{\it Numer.\ Math.\ Theor.\ Meth.\ Applic.\ }}
\def\PAMS{{\it Proc.\ Amer.\ Math.\ Soc.\ }}
\def\SIAMJNA{{\it SIAM J.\ Numer.\ Anal.\ }}
\def\SIAMR{{\it SIAM Rev.\ }}
\def\SIVP{{\it Signal, Image Video Proc.\ }}
\def\SR{{\it Soft Robotics }}

\begin{flushleft}

\end{flushleft}
\end{document}

%% file: section01.tex
\section{Introduction}

Planar and spatial polynomial Pythagorean-hodograph (PH) curves are characterized by a \emph{polynomial} parametric speed, and can be effectively represented in terms of complex and quaternion algebra, respectively \cite {farouki08}. We refer to \cite{review19} for a comprehensive recent introduction to PH curves, as well as to related constructions and applications. The spline extension of PH structures offers many advantages in the context of path generation, as for examples the possibility of easily constructing flexible fair shapes, and a piecewise polynomial arc length, a distinctive feature for motion control on a given path, see e.g., \cite{giannelli16} and references therein. 

The focus of this paper is on the construction of a \emph{local} interpolation scheme based on spatial PH quintic splines for given $C^1$ Hermite (positions and first-order derivatives) \emph{data streams}. The availability of input Hermite data is usually assumed when higher shape control for suitable path identification should be enabled. More specifically, each spline segment interpolates second and first order Hermite data at the initial and final end-point, respectively. In the spline extension of the scheme a $C^2$ smooth connection between successive spline segments is obtained by taking the locally required second-order derivative information from the previous segment. Consequently, the data stream spline interpolant is globally $C^2$ continuous and can be constructed for arbitrary $C^1$ Hermite data configurations.


Interpolation of spatial Hermite data by PH splines is a challenging problem involving a \emph{family} of solutions whose shape is strongly influenced by the choice of certain free parameters.  While general $C^1$ Hermite data can always be interpolated by cubic polynomial splines, the PH condition reduces the available degrees of freedom and PH cubic interpolants with $C^1$ or $G^1$ continuity not always exist \cite{juettler99,kwon10,pelosi05}. To avoid this problem the PH spline degree can be raised to 5 \cite{farouki02b,farouki08c,sestini13} or the cubic spline segment can be subdivided in two (or more) segments \cite{bastl14a,sestini14}. In a similar way, spatial $C^2$ PH spline interpolants can be (always) obtained either with PH curves of higher (nine) degree \cite{sir07} or with low-degree PH spline arcs of different kind. For example, the interpolation scheme proposed in \cite{bastl14b} relies on PH quintic triarcs, since PH quintic biarcs are not flexible enough for addressing the symmetric second-order Hermite interpolation problem. 

An increasingly range of applications nowadays requires a real-time processing of input data streams, usually consisting of position information with or without first  and second order Hermite data. The access to the whole data set is then not originally available and the streaming algorithm should be completely local and suitably exploit first (and second) derivative estimates to properly compute smooth paths with different order of continuity. As a consequence, symmetric interpolation algorithms which generate the same geometric path when the order of the whole data sequence is reverted (a desirable feature when the data are all simultaneously available) are not strictly needed in this context, see for example \cite{Debski1}.

Our method relies on PH quintic biarcs, since they ensure enough flexibility for addressing the considered interpolation problem, while preserving the practicality of PH quintics. Between any pair of successive Hermite data the corresponding path segment is defined as a spatial PH quintic biarc, which also requires a second derivative information at one of the two extrema. However, in the spline implementation of the scheme for Hermite data streams, this information is not required as additional input data since it is directly obtained by evaluating  the second derivative of the previous biarc at the joint point to guarantee curvature continuity at this point. A single PH quintic arc interpolating only first-order Hermite data is considered at the beginning to completely avoid the need of second-order derivative information. A simple and effective selection of the free parameters that arise in the interpolation problem is proposed by relying on the results for the $C^1$ PH quintic spline interpolants proposed in \cite{farouki08c}. The theoretical analysis proves the fourth approximation order of the local scheme when exact first order Hermite data and a second order error on the $C^2$ condition are considered. A selection of numerical examples confirms the same approximation for the spline extension of the scheme. In addition, the performance of the algorithm are also validated by considering an application-oriented point stream with automatically generated first-order derivative information.

The plan for this paper is as follows. The data stream interpolation algorithm based on polynomial Pythagorean-hodograph biarcs of degree five is presented in Section~\ref{sec:PH}. An effective strategy for the selection of the free parameters is proposed in Section~\ref{sec:free}. The approximation order of the local scheme is studied in Section~\ref{sec:approx}, while Section~\ref{sec:numerics} presents a selection of computed examples. Finally, Section~\ref{sec:closure} summarizes the key results of this paper.

%% file: section02.tex
\section{The interpolation scheme}
\label{sec:PH}

In this section we introduce the local interpolation problem addressed in the paper, referring to the Appendix for a short introduction to the quaternion algebra. We aim to construct a parametric PH quintic spline curve $\x(u), u\in[u_i\,,\,u_f]$ interpolating two assigned points $\p_i\,\,,\, \p_f \in \E^3$ at the end parameter values and such that
\begin{equation} \label{asymprobder}
\left\{
\begin{array}{llll}
\x'(u_i) &= \vv_i\,, &  \x'(u_f) &= \vv_f\,, \cr
\x''(u_i) &= \w_i\,, & \ \cr 
\end{array}
\right.
\end{equation}
where $\vv_i, \vv_f$ and $\w_i$ are assigned vectors\footnote{
Note that in the spline interpolation of an Hermite data stream, the vector $\w_i\,,$ is not an additional input information since it is simply taken from the previous spline segment of the path.}  in $\R^3$ and the $'$ symbol denotes derivatives with respect to the global parameter $u.$  In particular, we consider a PH quintic biarc, a spline curve composed by two PH quintic polynomial segments $\x_i,\, \x_f$ joining at the inner parameter value $u_m,\, u_i < u_m < u_f,$
$$\x(u) = \left\{\begin{array}{ll} 
\x_i(u) & \mbox{for } u \in [u_i\,,\,u_m], \cr
\x_f(u) & \mbox{for } u \in [u_m\,,\,u_f]. \cr
\end{array} \right.$$
Let $t = (u-u_i)/(u_f - u_i)$ be the local parameter varying in $[0\,,\,1]$ associated to the whole biarc, and $\hat t = (u_m-u_i)/(u_f-u_i)$ be the parameter value where the two PH quintic segments of the biarc join, i.e.
$$\hat t := \frac{h_i}{h_i + h_f}\,,$$
with $h_i := u_m-u_i$ and $h_f := u_f-u_m.$ 
In addition, let $\tau$ and $\eta$ denote the local parameter values both varying in $[0\,,\,1]$ associated with $\x_i$ and $\x_f,$ respectively,
$$\tau := \frac{u-u_i}{h_i} = \frac{t}{\hat t}\,, \qquad \eta := \frac{u-u_m}{h_f} = \frac{t-\hat t }{1 - \hat t}\,.$$
In order to consider PH quintic biarcs, we set
\begin{equation} \label{hodo}
\frac{d \x_i}{d \tau}(\tau) := \A(\tau)\, \i \,\A^*(\tau)\,, \qquad \frac{d \x_f}{d\eta}(\eta)  := \B(\eta)\, \i \,\B^*(\eta)\,,
\end{equation}
where 
\begin{equation} \label{preimage}
\A(\tau) := \sum_{j=0}^2 \A_j B_j^2(\tau)\,, \qquad \B(\eta)  := \sum_{j=0}^2 \B_j B_j^2(\eta)\,,
\end{equation}
are quadratic quaternion polynomials which define the {\it pre-image} of $\x_i$ and $\x_f,$ in the Bernstein basis, with $B_0^2(\xi) := (1-\xi)^2, B_1^2(\xi) := 2\xi(1-\xi)$ and $B_2^2(\xi) := \xi^2$ and $\A_j\,, \B_j, j=0,1,2$ denoting quaternion coefficients belonging to $\mathbb{H}.$ As well known, the two quintic arcs $\x_i$ and $\x_f$ can be written as Bezier curves as follows,
\begin{equation} \label{PH}
\x_i(\tau) =  \sum_{j=0}^5 \q_{ji} B_j^5(\tau)\,, \qquad \x_f(\eta) =  \sum_{j=0}^5 \q_{jf} B_j^5(\eta)\,,
\end{equation}
where $\q_{0i} = \p_i,$ $\q_{5f} = \p_f,$ and
\begin{equation} \label{cpi}
\begin{array}{ll}
\q_{1i} &= \q_{0i} + \frac{1}{5} 
\A_0\,\i\,\A_0^*, \cr 
\q_{2i}  &= \q_{1i} + \frac{1}{10}
\left(\A_0\,\i \,\A_1^* + 
\A_1\,\i \,A_0^*\right),  \cr
\q_{3i} & = \q_{2i} + \frac{1}{30}
\left(\A_0\,\i\,\A_2^*+ 4\,\A_1\,\i\,\A_1^*+ \A_2\,\i\, \A_0^*\right),\cr
\q_{4i} & =\q_{3i} + \frac{1}{10}
\left(\A_1\,\i\,\A_2^*+ \A_2\,\i\, \A_1^* \right),\cr 
\q_{5i}  &= \q_{4i} + \frac{1}{5} \A_2\,\i\,\A_2^*\,,
\end{array} 
\end{equation}
and
\begin{equation} \label{cpf}
\begin{array}{ll}
\q_{4f} &= \q_{5f} - \frac{1}{5} 
\B_2\,\i\,\B_2^*,\cr
\q_{3f}  &= \q_{4f} - \frac{1}{10} \left(\B_1\,\i \,\B_2^* + \B_2\,\i \,\B_1^* \right), \cr
\q_{2f} & = \q_{3f} - \frac{1}{30}
\left( \B_0\,\i\,\B_2^*+ 4\,\B_1\,\i\,\B_1^*+ \B_2\,\i\, \B_0^* \right),\cr
\q_{1f} & =\q_{2f} - \frac{1}{10} \left(\B_0\,\i\,\B_1^*+ \B_1\,\i\, \B_0^*\right),\cr
\q_{0f}  &= \q_{1f} - \frac{1}{5}
\B_0\,\i\,\B_0^*\,. 
\end{array}
\end{equation}
Since the position interpolation conditions are ensured by the assumption $\q_{0i}=\p_i$ and $\q_{5f} = \p_f,$ it is clear that we need to define the  six quaternion coefficients ${\cal A}_j, \, {\cal B}_j,$ $ j=0,1,2$ in order to satisfy the conditions in (\ref{asymprobder}). Note that the scalar degrees of freedom available for compute the solution are just $23$ (and not $6 \cdot 4 =24$) since a spatial PH curve does not change if any quaternion coefficient of its pre-image is right multiplied by a common complex unit factor of the form $\cos \theta + \i \sin \theta.$ Being the considered smoothness and interpolation conditions just $6 \cdot 3 =18$ (two vector conditions on first and second derivatives at the left end point, one vector condition on first derivative at the right end point  and three vector conditions to ensure the $C^2$ joint at $u_m$ of the biarc), there are necessarily some free parameters in the scheme. Additional conditions to identify a suitable $C^2$ interpolating biarc between the family of formal solutions should then be considered. This has allowed us a preliminary  removal of two free parameters by requiring that the scheme produces just one PH quintic instead of a biarc, whenever possible. We then impose a $C^1$ joint between the quaternion pre-images of the two polynomial segments of the biarc. This also ensures $C^1$ smoothness at the same point to the Euler--Rodrigues frame associate to the PH biarc, see for example \cite{review19} and references therein. If, for simplicity, $u_m$ is chosen at the midpoint between $u_i$ and $u_f,$ this condition corresponds to assign ${\cal B}_0$ and ${\cal B}_1$ as follows,
 \begin{equation} \label{suff}
  \B_0 =   \A_2\,, \qquad \B_1 = 2 \A_2 - \A_1 \,.
  \end{equation}
  
\begin{figure}[!t]
\centering

\subfigure[$a_1 = -2.0, -1.5 , \ldots ,2.0 , 2.5. $ ]{
\includegraphics[width=0.45\linewidth]{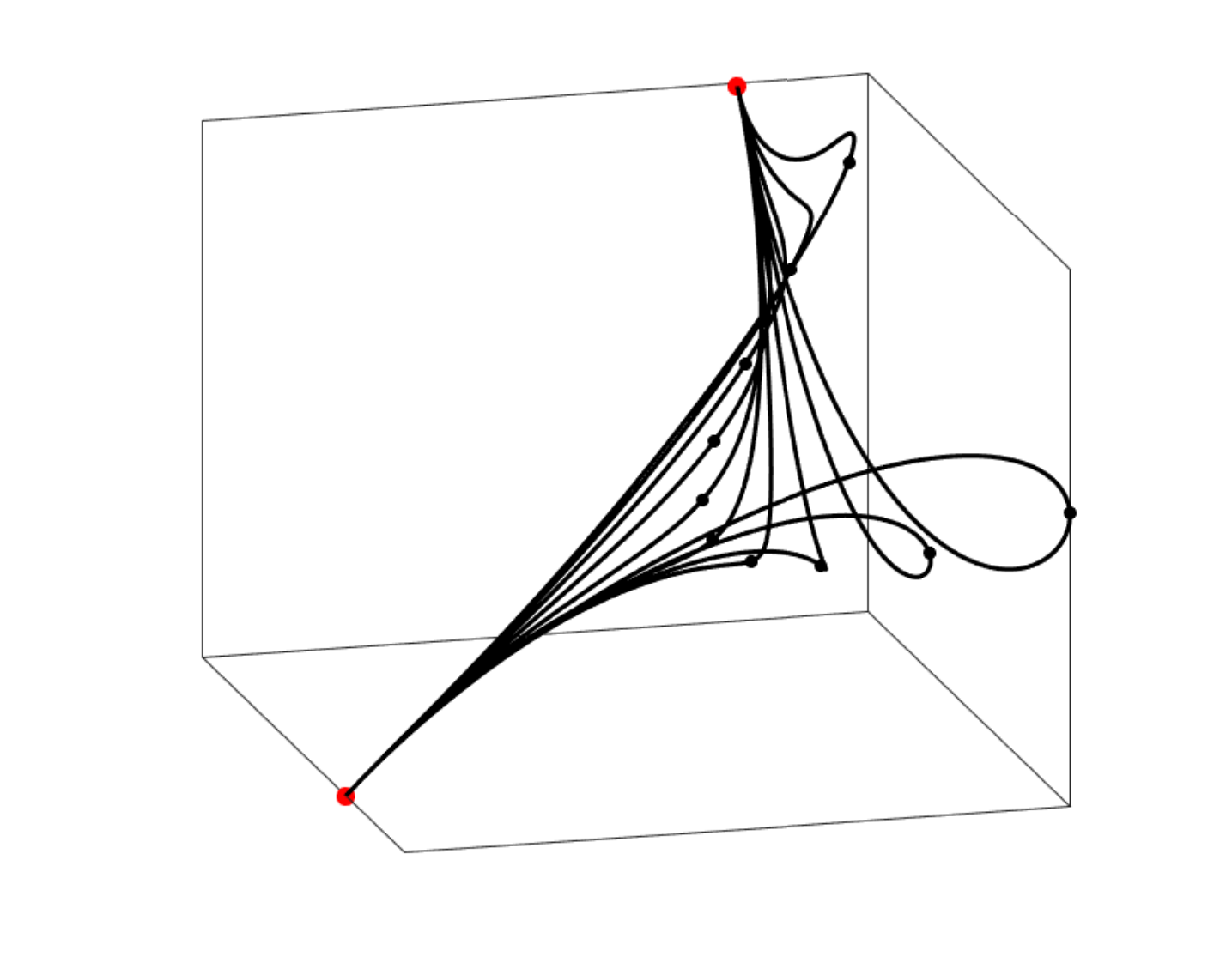}}
\subfigure[${\alpha}_2 = 0, \frac{\pi}{5}, \ldots , \frac{8 \pi}{5},\frac{9 \pi}{5}.$]{
\includegraphics[width=0.45\linewidth]{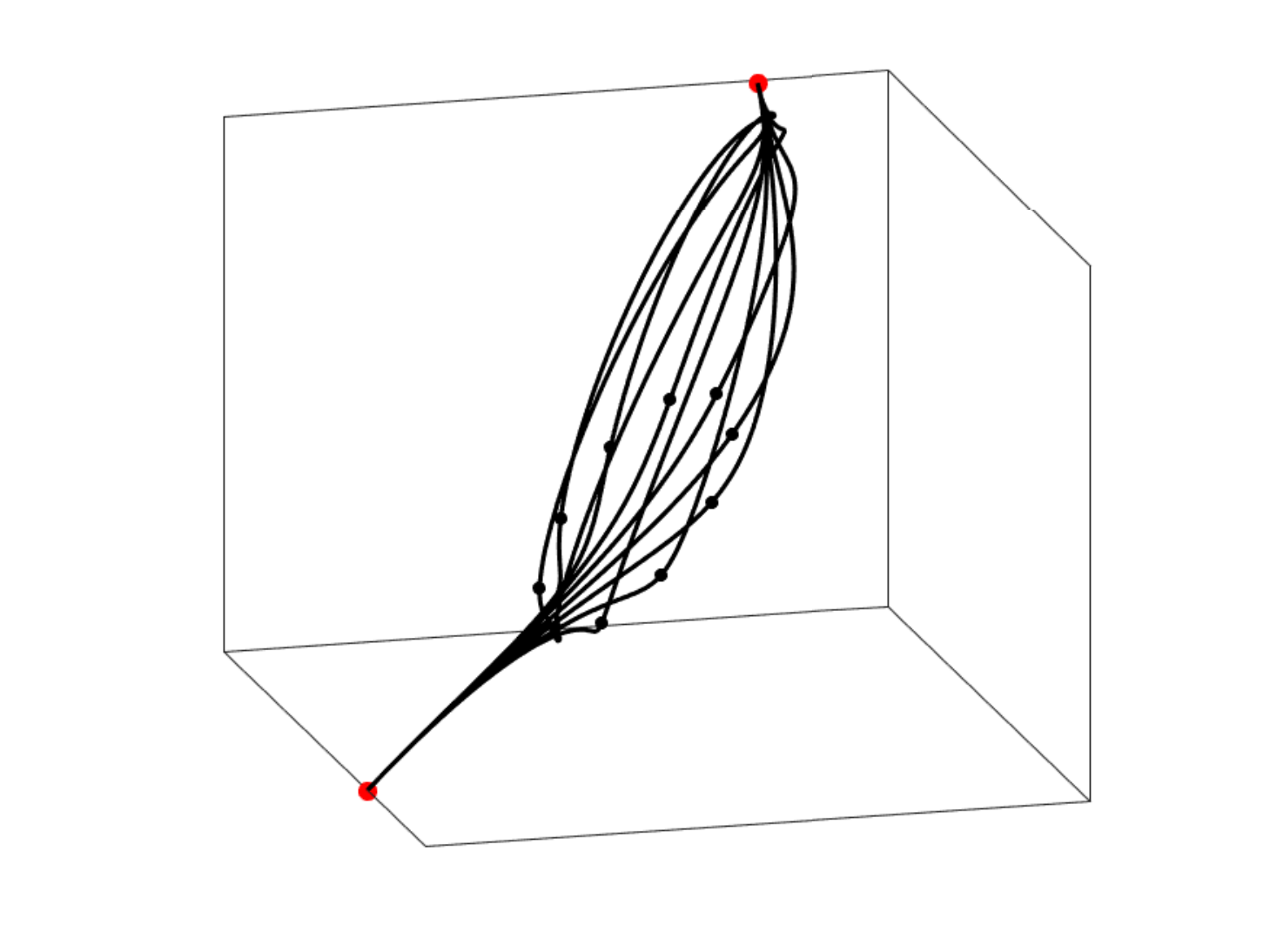}}

\caption{Influence of the parameters $a_1$ and $\alpha_2$ on the shape of the biarc. The value of the other parameters is fixed to the one resulting from the algorithm that will be presented below. The Hermite conditions are: $\p_i = (0 , 0 , 0)^T $ , $\p_f = (1 , 1 , 1)^T $, $ \vv_i = (1 , 0 , 1)^T $ , $ \vv_f = (0 , 1 , 1)^T $, $ \w_i = (-0.1 , 0.5 , -1.5)^T $.}
\label{fig:param_var}
\end{figure}

The derivative interpolation requirements assigned at the end points are the following,
 \begin{equation} \label{interp} \frac{d\x_i}{d u}(u_i) = \vv_i\,,\quad  \frac{d\x_f}{d u}(u_f) = \vv_f\,, \quad  \frac{d^2\x_i}{d u^2}(u_i) = \w_i.
 \end{equation}
 Now, considering $\hat t = 0.5\,,$ that is $h_i = h_f = h \,:=\, (u_f-u_i)/2\,,$ the following two derivative chain rules can be easily obtained,
 $$\frac{d\x_i}{d u} = \frac{1}{2h}\frac{d\x_i}{d t} = \frac{1}{h}\frac{d\x_i}{d \tau}\,, \quad \frac{d\x_f}{d u} = \frac{1}{2h}\frac{d\x_f}{d t} =  \frac{1}{h}\frac{d\x_f}{d \eta}\,,$$  
 $$\frac{d^2\x_i}{d u^2} = \frac{1}{4h^2}\frac{d^2\x_i}{d t^2} = \frac{1}{h^2}\frac{d^2\x_i}{d \tau^2}\,, \quad \frac{d^2\x_f}{d u^2} = \frac{1}{4h^2}\frac{d^2\x_f}{d t^2} =  \frac{1}{h^2}\frac{d^2\x_f}{d \eta^2}\,.$$
Consequently, using the derivative formula for polynomials in Bernstein form, the conditions in (\ref{interp}) correspond to the following vector conditions in the quaternion algebra,
 \begin{equation} \label{endcond}
 \left\{ \begin{array}{l} \A_0\,\i\,\A_0^* = h \vv_i\,, \cr
 \B_2\,\i\,\B_2^* = h  \vv_f\,, \cr
 2 \left( (\A_1 - \A_0) \,\i\,\A_0^* +A_0 \,\i\,(\A_1 - \A_0)^*\right) = h^2 \w_i\,\, \end{array} 
   \right.
  \end{equation}   
where, taking into account the first two equations, the last one reduces to
 \begin{equation} 
 \label{deriv2}
 2(\A_1\,\i\,\A_0^* + \A_0\,\i\,\A_1^*) = h_i^2 \w_i + 4h_i \vv_i\,.
 \end{equation} 
From the general quaternion solution of a quadratic equation \cite{farouki02b}, we derive the expressions 
 \begin{equation} \label{A0B2}
  \A_0 = \sqrt{h_i \vert \vv_i \vert} \ \U_0  \,, \quad \B_2 = \sqrt{h_f \vert \vv_f \vert}\  \U_2\,,
  \end{equation}
for $\A_0$ and $\B_2$ which satisfy the first two conditions in (\ref{endcond}), with the unit quaternions
 \begin{equation} \label{U0U2} \U_0 := \u_i (\cos \alpha_0 + \i \sin \alpha_0) \,, \qquad \U_2 := \u_f (\cos \beta_2 + \i \sin \beta_2)\,,   \end{equation}
depending on the two free angular parameters $\alpha_0$ and $\beta_2$, and the unit vector quaternions
 $$ \u_i := \ \frac{\delta_i + \i}{\vert \delta_i + \i \vert}\,, \quad  \u_f := \ \frac{\delta_f + \i}{\vert \delta_f + \i \vert}\,,$$
written in terms of $\delta_i := \vv_i/\vert \vv_i \vert$ and $\delta_f := \vv_f/\vert \vv_f \vert\,.$
By setting 
 \begin{equation} \label{dfdef}\d_i :=  \frac{h_i}{4}\ \left(h_i \w_i + 4\vv_i \right),
 \end{equation}  
 the vector condition in (\ref{deriv2}) is fulfilled if 
 \begin{equation} \label{A1} \A_1 =  -\frac{(a_1 + \d_i )  \A_0 \i}{\vert \A_0 \vert^2} \,,
 \end{equation}
 where $a_1 \in \R$ is a scalar free parameter.  
 Now, considering (\ref{suff}), in order to ensure a $C^2$ connection at the joint point between the two polynomial segments of the  biarc we need only to require $\q_{0f} = \q_{5i},$  that is
 \begin{align} \label{lasteq}
& \p_i + \frac{1}{5} h \vv_i + \frac{1}{20} \left(h^2 \w_i + 4h \vv_i \right) + 
 \frac{1}{30} \left(\A_0\,\i\,\A_2^*+ 4\,\A_1\,\i\,\A_1^*+ \A_2\,\i\, A_0^*\right) + \frac{1}{10} \left(\A_1\,\i\,\A_2^*+ \A_2\,\i\, A_1^*\right)  \cr
& \qquad + 
\frac{1}{5} \A_2\,\i\,\A_2^*  \,=\, 
\p_f - \frac{1}{5} h \vv_f - \frac{1}{10}\left(\B_1\,\i \,\B_2^* + \B_2\,\i \,B_1^*\right) - 
 \frac{1}{30} \left(\B_0\,\i\,\B_2^*+4\,\B_1\,\i\,\B_1^*+ \B_2\,\i\, B_0^*\right) \cr
& \qquad - \frac{1}{10} \left(\B_0\,\i\,\B_1^*+ \B_1\,\i\, B_0^*\right)  - \frac{1}{5} \B_0\,\i\,\B_0^*\,.
\end{align}
After replacing the expressions for ${\cal B}_0$ and ${\cal B}_1$ given in (\ref{suff}), moving on the left hand side all the terms containing $\A_2$ and on the right all the others, this equation can be properly simplified to
\begin{align}
&\A_2\,\i\,\A_2^* + \frac{1}{40} \left(\G\,\i\,\A_2^* + \A_2\,\i\,\G^*\right) = \frac{15}{20} \c - \frac{2}{10} \A_1\,\i\,\A_1^* + \frac{3}{40} \left(\A_1\,\i\,\B_2^* + \B_2\,\i\,\A_1^*\right)\,,
\end{align} 
where 
\begin{equation} \label{Gdef} 
\G := \A_0-8\A_1+7\B_2 \,
\end{equation}
and
$$\c := (\p_f - \p_i) - \frac{1}{5} h (\vv_f + \vv_i) - \frac{1}{20} \left(h^2 \w_i + 4h \vv_i \right)\,$$
This last vector equation can also be further algebraically manipulated to arrive at the following final form,
\begin{equation} \label{eqfinal}
 \left(\A_2 +  \frac{1}{40} \G\right) \i \left(\A_2 +  \frac{1}{40} \G\right)^* = \b\,,
 \end{equation}
where
$$ \b : = \frac{15}{20} \c - \frac{2}{10  } \A_1\,\i\,\A_1^* + \frac{3}{40} \left(\A_1\,\i\,\B_2^* + \B_2\,\i\,\A_1^*\right) \,+\,  \frac{1}{1600} \G\, \i\, \G^*\,.
$$
Thus, considering again the general solution of a quadratic quaternion equation of type ${\cal V}\,\i\,{\cal V}^* = \r,$ from (\ref{eqfinal}) we get
\begin{equation} \label{A2} 
\A_2 \,=\, - \frac{1}{40} \G + \q \ (\cos \alpha_2 + \i \sin \alpha_2)\,,
\end{equation}
where $\alpha_2$ is another free angular parameter and for brevity we have set,
\begin{equation}
\label{qdef}
 \q :=\sqrt{ \vert \b \vert} \ \frac{\i + \frac{\b}{\vert \b \vert}}{\left\vert \i + \frac{\b}{\vert \b \vert} \right\vert}\,.
 \end{equation}
Summarizing, the scheme has four free parameters: the three angles $\alpha_0, \beta_2,$ and $\alpha_2$, together with the real coefficient $a_1.$ However, the next proposition proves that the real shape angular parameters are just two: $\alpha_2 - \alpha_0$ and $\beta_2-\alpha_0,$
\begin{prn}
For any $\theta \in [0\,,\,2\pi)$ the replacement of $\alpha_0, \alpha_2$ and $\beta_2$ in our scheme with $\alpha_0 +\theta, \alpha_2 +\theta$ and $\beta_2+\theta$, respectively, produces the same PH quintic biarc.
\end{prn}
\prf We need to prove that adding $\theta$  to the free angular parameters $\alpha_0, \alpha_2$ and $\beta_2$ implies that all the quaternion coefficients ${\cal A}_j, \, {\cal B}_j,$ $ j=0,1,2$ are right multiplied by the complex unit quaternion $e^{\i\, \theta} = \cos \theta + \i \sin \theta.$ For ${\cal A}_0$ and ${\cal B}_2,$ this can be verified using equations (\ref{A0B2}) and (\ref{U0U2}) together with trigonometric sum formulas.  As for ${\cal A}_1,$ considering that we have already proved that ${\cal A}_0$ is right multiplied for $e^{\i\,\theta},$ recalling that $e^{\i\,\theta}\i = \i\, e^{\i\,\theta},$ formula (\ref{A1}) implies the result. As a consequence, formula (\ref{Gdef}) implies that also the quaternion ${\cal G}$ is right multiplied by this factor. Consequently, formulas (\ref{A2}) and (\ref{suff}) imply that this is also true for ${\cal A}_2$ and ${\cal B}_0, {\cal B}_1$, respectively. \QED \\
As usual with PH curves, the shape parameters highly influence the  resulting interpolant, as shown in Figure \ref{fig:param_var}. Reasonable and effective criteria for their choice are necessary. More precisely, since the scheme has to be used for defining a spline curve, it is fundamental that just one PH biarc is locally identified. This means that we need to properly fix our three shape parameters, $\alpha_2 - \alpha_0, \beta_2-\alpha_0$ and $a_1.$ In addition, since we are interested in defining an interpolation scheme suitable for real-time applications, it is also necessary that these criteria have an easy implementation. In the next section we introduce a simple and effective selection strategy driven by a reference PH quintic curve  introduced in \cite{farouki08c}  to solve the $C^1$ Hermite interpolation problem. Successively, in Section~\ref{sec:free}, we show that this strategy is reasonable from the approximation point of view, since it ensures fourth approximation order to our scheme. This is the same approximation order characterizing the $C^1$ PH spline quintic Hermite scheme used as a reference for the free parameter selection.

%% file: section03.tex
\section{Selection of free parameters}
\label{sec:free}
As already mentioned in the previous section, the free angular parameters $\alpha_2-\alpha_0, \beta_2- \alpha_0$ and the real parameter $a_1$ of our approach can highly affect the shape of the interpolating biarc.   In this section we introduce our data-dependent selection strategy  based on a reference PH curve easy to be constructed. This curve is defined as a PH quintic interpolating the given zero and first order data at the extrema but not the assigned second order information at the left end point (this would not be possible in the general case, since a PH quintic is not flexible enough for ensuring this additional condition). Now, as shown in \cite{farouki08c}, the determination of an Hermite PH quintic interpolant is affected by two free shape parameters (all of angular type)  and four data-dependent different criteria for their selection were there introduced and compared. Here we rely on the so called CC criterion\footnote{The CC label is the acronym for Cubic-Cubic, since the selection strategy of both the angular shape parameters affecting an Hermite PH quintic interpolant is done in this case using the standard cubic Hermite interpolant as a reference curve.} whose implementation is very easy. Besides having a very good behavior for general non asymptotic data, in \cite{sestini13} it has also been proved that this criterion ensures fourth approximation order when the  scheme is applied to the reconstruction of a given analytic smooth curve. 

\subsection{The reference PH quintic} \label{reference}
Let us denote with
$${\cal V}^H(t) = \sum_{j=0}^2 {\cal V}_j^H B_j^2(t)\,, \qquad t \in [0\,,\,1]\,,$$
the pre-image expressed in the local parameter $t$ of the CC PH quintic polynomial segment $\x^H(t)$ verifying the assigned end Hermite conditions,
$$ \x^H(0) = \p_i \,, \quad \frac{d\x_i}{dt}(0) = 2h \vv_i \,, \quad \x^H(1) = \p_f \,, \quad \frac{d\x_f}{dt}(1) = 2h \vv_f \,.$$
Concerning the algorithm for its construction, we recall some basic points, referring to \cite{farouki08c} for the details. First we observe that, dealing with PH curves, one of the two end point interpolation conditions come for free by taking into account the integral definition of a PH curve from its hodograph. Consequently, only three of the four Hermite conditions are available for the characterization of the three quaternion coefficients ${\cal V}_j^H, j=0,1,2$ and a free angular parameter is associated to each of them. On the other hand, only two of these three free angles are real shape parameters, since the interpolating PH quintic does not change if all the quaternion coefficients are multiplied by a trigonometric factor of the form $\cos\theta + \i \sin \theta.$ For all the criteria proposed in \cite{farouki08c}, the angular free parameter associated to ${\cal V}_1^H$ was chosen equal to zero and the other two were selected according to the considered criterion. We refer to \cite{farouki08c} for the details on the CC criterion by simply outlining that it allows degree reduction of the pre-image and, consequently, the definition of a PH cubic Hermite interpolant, whenever possible. 

Since we use the CC PH quintic Hermite interpolant as a reference curve to construct our PH quintic biarc, it is useful to represent it as a biarc with joint parameter at $t = \hat t = 0.5.$ We then set 
\begin{equation}
\x^H(t) = \left\{\begin{array}{ll} 
\x_i^H(t) & \mbox{for } t \in [0\,,\,\hat t], \cr
\x_f^H(t) & \mbox{for } t \in [\hat t\,,\,1]. \cr
\end{array} \right.
\label{eq:PHsub}
\end{equation}
Using again on the left the local parameter $\tau \in [0\,,\,1]$ and on the right $\eta \in [0\,,\,1]$ we observe that 
$$\frac{d\x_i^H}{d\tau} =  {\cal A}^H(\tau)\, \i\, {\cal A}^{H*}(\tau)\,, \quad \frac{d\x_f^H}{d\eta} =  {\cal B}^H(\eta) \,\i \,{\cal B}^{H*}(\eta)\,,$$
where
$${\cal A}^H(\tau) = \sum_{j=0}^2 {\cal A}_j^H B_j^2(\tau) \,, \quad  {\cal B}^H(\eta) = \sum_{j=0}^2 {\cal B}_j^H B_j^2(\eta) \,, $$
with
\begin{align} \label{quatref} 
&{\cal A}_0^H \,:=\, \frac{1}{\sqrt{2} } {\cal V}_0^H\,, \,\, {\cal A}_1^H \,:=\,   \frac{1}{2\sqrt{2} }( {\cal V}_0^H + {\cal V}_1^H)\,,\,\,
{\cal B}_1^H \,:=\,   \frac{1}{2\sqrt{2} }( {\cal V}_1^H + {\cal V}_2^H)\,, \,\, {\cal B}_2^H \,:=\, \frac{1}{\sqrt{2} } {\cal V}_2^H\,,   \cr
\ \cr
&{\cal A}_2^H = {\cal B}_0^H := \frac{1}{4\sqrt{2}} ( {\cal V}_0^H + 2 {\cal V}_1^H + {\cal V}_2^H)\,.
\end{align}

\subsection{The parameter selection strategy for the biarc}

The quaternion coefficients ${\cal A}_j^H, \, {\cal B}_j^H, j=0,1,2$ introduced  in (\ref{quatref}) to express the hodograph of the reference CC PH quintic Hermite interpolant in biarc form are directly used in our strategy to drive the selection of the free parameters. We observe that, even if only the differences $\alpha_2-\alpha_0$ and $\beta_2-\alpha_0$ are real angular shape parameters, in the algorithm we specify all the three angles $\alpha_0, \beta_2, \alpha_2,$ since we refer to a specific pre-image of the reference PH curve for their selection. Note also that first $\alpha_0$ and $\beta_2$ are selected referring to ${\cal A}_0^H$ and ${\cal B}_2^H,$ respectively. Subsequently, the real parameter $a_1$ involved in the definition of ${\cal A}_1$ is determined through an explicit formula which only involves ${\cal A}_0$ and ${\cal A}_1^H.$ Finally, the remaining free angle $\alpha_2$ is determined with an explicit formula depending on all the previous choices, as well as on the reference quaternion coefficient ${\cal A}_2^H$.

The choice of the two extreme angular free parameters $\alpha_0$ and $\beta_2$ is very easy, since we just set
$${\cal A}_0 = {\cal A}_0^H\,, \qquad {\cal B}_2 = {\cal B}_2^H\,.$$
We then determine the free parameter $a_1$ in order to minimize the following scalar quantity,
$$\vert \A_1 - \A_1^H \vert^2 = \A_1\A_1^* - \left(\A_1 \A_1^{H*} + \A_1^H\A_1^*\right) + \A_1^H \A_1^{H*}\,.$$
Recalling that $\A_1 = (a_1 + \d_i)\, \U_0/\vert \A_0 \vert\,,$ with $\U_0 = -(\A_0 \i)/\vert\A_0\vert,$ we get
$$\vert {\cal A}_1 - {\cal A}_1^H \vert^2 = \frac{a^2_1}{\vert \A_0 \vert^2}  - \frac{a_1}{\vert \A_0 \vert} \left(\U_0\A_1^{H*} + \A_1^H \U^*_0 \right) + K$$
where $$K := \frac{\vert \d_i\vert^2}{\vert \A_0 \vert} -\frac{\left(\d_i \U_0 \A_1^{H*} - \A_1^H \U_0^*\d_i  \right)}{\vert \A_0 \vert} + \A_1^H \A_1^{H*}\,. $$
Thus the minimum is obtained by setting
 \begin{equation} \label{a1def}
 a_1 = \frac{1}{2}\vert \A_0 \vert \left(\U_0 \A_1^{H*} + \A_1^H  \U^*_0 \right) \,=\, \frac{1}{2}  \left(\A_0 \A_1^{H*} + \A_1^H  \A^*_0 \right)\,.
 \end{equation}
Finally, we determine the last free parameter $\alpha_2$ by minimizing the following quantity,
\begin{equation}\label{eq:min}
\vert \A_2 - \A_2^H \vert^2 = \A_2 \A_2^* - \left(\ A_2 \A_2^{H*} + \A_2^H\A_2^*\right) + \A_2^H \A_2^{H*}\,.
\end{equation}
Since $\A_2$ depends on $\alpha_2$ and, in particular, 
$$\A_2 = -\frac{1}{4(10 + b_1^2)} \G +\q e^{\i\alpha_2}\,,$$ we have
\begin{align*}
 \A_2\A_2^* & = \left( -\frac{1}{4(10 + b_1^2)} \G +\q\, e^{\i\,\alpha_2}\right)\left( -\frac{1}{4(10 + b_1^2)}\, \G^* - e^{-\i\,\alpha_2}\q\right) \cr
 & = \frac{1}{16(10 + b_1^2)^2}\, \G\G^* -\q^2 
 + \frac{1}{4(10 + b_1^2)} \left(\G e^{-\i\,\alpha_2}\q -\q e^{\i\,\alpha_2}\G^* \right)
\end{align*}
and
\begin{align*}
\left(\A_2 \A_2^{H*} + \A_2^H\A_2^*\right) 
=  -\frac{1}{4(10 + b_1^2)} \left(\G \A_2^{H*} + \A_2^H\G^*\right) +
 \left(\q e^{\i\, \alpha_2}\A_2^{H*} - \A_2^H e^{-\i\,\alpha_2} \q \right) 
\end{align*}
This implies that
$$\vert \A_2 - \A_2^H \vert^2 = f_1 \cos \alpha_2 + f_2 \sin \alpha_2 + g  \,,$$
with 
$$g := \frac{1}{16(10 + b_1^2)^2} \G\G^* -\q^2  +\frac{1}{4(10 + b_1^2)} \left(\G \A_2^{H*} + \A_2^H\G^*\right) + \A_2^H \A_2^{H*}$$ 
and
\begin{align} \label{fcoefdef}
f_1 &:= \frac{1}{4(10 + b_1^2)} (\G  \q -\q  \G^* ) - \left(\q \A_2^{H*} - A_2^{H} \q \right), \cr
f_2 &:= -\frac{1}{4(10 + b_1^2)} (\G \i \q +\q \i \G^* ) - \left(\q \i \A_2^{H*} + \A_2^H \i \q\right)\,.
\end{align} 
Thus we minimize \eqref{eq:min} if
\begin{equation} \label{alpha2def}
\alpha_2 = \pi + \atantwo(f_2,f_1).
\end{equation}

\begin{figure}[!t]
\centering

\subfigure[$a_1$ variation plot.]{
\includegraphics[width=0.45\linewidth]{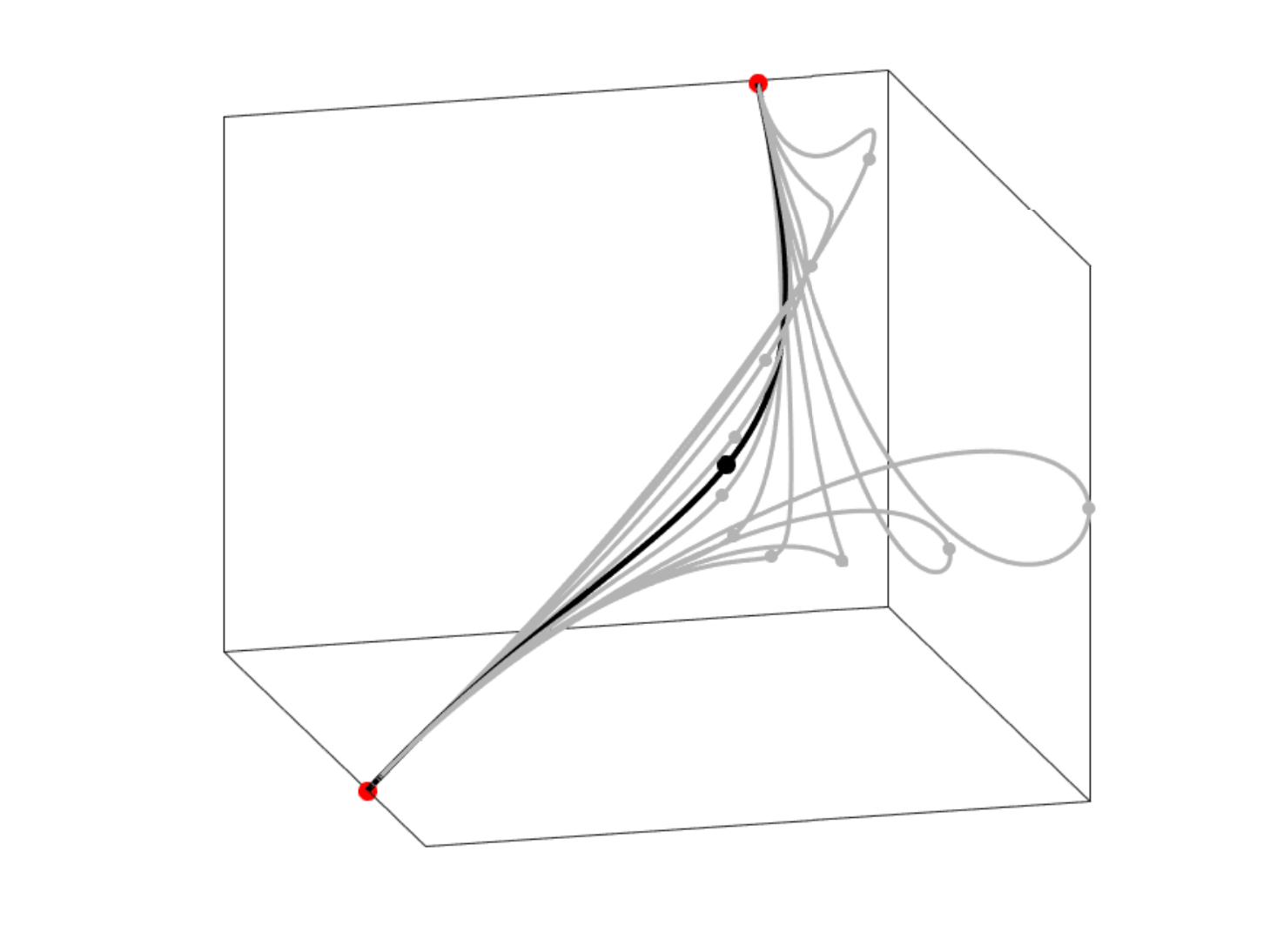}}
\subfigure[${\alpha}_2$ variation plot.]{
\includegraphics[width=0.45\linewidth]{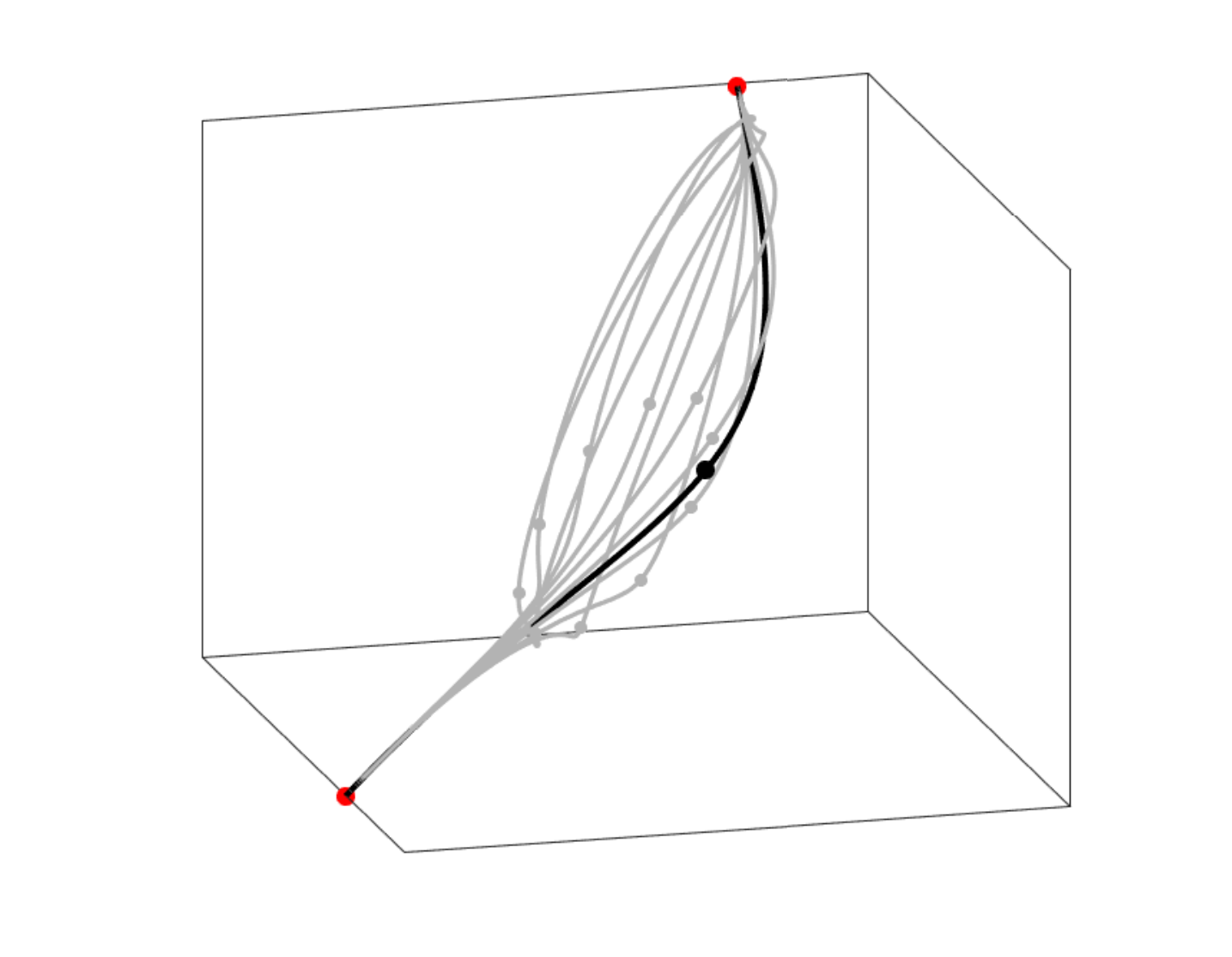}}

\caption{Influence of $a_1$ and $\alpha_2$ parameters on the shape of the biarc. The curve obtained with our selection strategy is shown in bold.}
\label{fig:param_var2}
\end{figure}

In Figure~\ref{fig:param_var2} the biarc resulting from this particular choice of the parameters is plotted together with the same set of variation of Figure~\ref{fig:param_var}.

%% file: section04.tex
\section{Approximation Order}
\label{sec:approx}

This section presents the asymptotic analysis to prove that the  $C^2$ PH quintic Hermite interpolation scheme here introduced has fourth approximation order. This is the same approximation order characterizing the $C^1$ first order PH quintic Hermite interpolation scheme introduced in \cite{farouki08c}, provided that the CC criterion is used to fix its two free angular parameters. This theoretical result was proved in \cite{sestini13} developing a local asymptotic analysis, i.e. just for one spline segment. The symmetric end-point first order Hermite data were there taken from an infinitesimal portion of a given smooth curve $\r(s)$ to be approximated, where $s$ is the arc-length parameter varying from $0$ to $\Delta s.$ Note that in \cite{sestini13} the fourth approximation order  can be trivially extended to the spline formulation of the scheme, since it computes each segment of the PH quintic spline independently from the others with the same local approach. 

We here present an analogous asymptotic analysis for our new $C^2$ forward scheme, relying on a MAPLE symbolic implementation of the local scheme. However, the situation in the spline formulation is now different: to obtain a $C^2$ interpolant, the biarc segments have to be computed in sequential order by using the additional second order information at the first end point of a segment from the previous one. In the local symbolic implementation of the scheme it is then reasonable to assume that this information can be affected by an error with an $O(\Delta s^2)$ expansion. Under this assumption we are able to prove that locally also the new scheme has fourth approximation order. A theoretical proof of the convergence order of the spline formulation of the scheme would require in this case also the analysis of the error propagation on the second derivative information at the joint points between successive biarcs. Even if this point was not formally addressed in virtue of the difficulties strictly connected to the non linear nature of any PH interpolation scheme, all the numerical experiments confirm the fourth approximation also to the spline formulation of the scheme.


Let us consider the portion of a given smooth curve $\r(s)$ restricted to $[0\,,\,\Delta s],$ and the symbolic computation of the PH biarc defined in Sections~\ref{sec:PH} and \ref{sec:free}  using the following input data:  
\begin{align}
\label{eq:HermiteConditions}
\nonumber
\p_i &= \r(0), &\p_f &= \r(\Delta s ), \\ 
\vv_i &= \frac{d\r}{ds}(0) , &\vv_f &= \frac{d\r}{ds}(\Delta s ), \\
\nonumber
\w_i &= \frac{d^2\r}{ds^2}(0) + O(\Delta s^2).
\end{align}

The following proposition is based on symbolic computations performed in the MAPLE computing environment (the MAPLE worksheet is available from the authors upon request).
\begin{prn}
Let us consider an arc-length parametrized curve $\r(s)$ with $s \in [0 , \Delta s]$ and $C^{10}$ continuity and let $\x$ be the PH quintic biarc fulfilling the vector interpolation conditions in (\ref{eq:HermiteConditions}). We have
\begin{equation*}
\Vert \x(t) - \r(t \Delta s)\Vert_2 = O(\Delta s^4), \,\,\,\, \forall t \in [0 , 1] \,\, \mbox{and }\,\, \Vert \frac{d^2\x}{ds^2}(\Delta s) - \frac{d^2\r}{ds^2}(\Delta s) \Vert_2 \,=\, O(\Delta s^2).
\end{equation*}
\end{prn}
\prf Denoting for brevity the Euclidean norm $\Vert \cdot \Vert_2$ just with $\Vert \cdot \Vert$ and using the triangular inequality, it is possible to write
\begin{equation*}
\Vert \x(t) - \r(t \Delta s)\Vert \leq \Vert \x(t) - \x^H(t)\Vert + \Vert \x^H(t) - \r(t \Delta s)\Vert
\end{equation*}
where $\x^H$ is the CC PH quintic interpolant introduced in \cite{farouki08c}. Now, under the same hypotheses here assumed, in \cite{sestini13} it has been proved that $\x^H$ provides a fourth order accurate approximation of $\r$. We can then study the first term of the previous inequality.
Considering the biarc representation of $\x^H$ introduced in (\ref{eq:PHsub}),  we can write
\begin{equation*}
\Vert \x(t) - \x^H(t)\Vert = \left\{\begin{array}{ll} 
\Vert \x_i(t) - \x_i^H(t) \Vert & \mbox{for } t \in [0\,,\,0.5], \cr
\Vert \x_f(t) - \x_f^H(t) \Vert & \mbox{for } t \in [0.5\,,\,1]. \cr
\end{array} \right. \end{equation*}
Since Bernstein polynomials are nonnegative in $[0\,,\,1]$ and sum up to 1, we can derive the following inequalities
$$\Vert \x_i(t) - \x_i^H(t) \Vert \le  \max_{k=0,\ldots,5} \Vert \q_{ik} - \q_{ik}^H \Vert\,, \quad
\Vert \x_f(t) - \x_f^H(t) \Vert \le  \max_{k=0,\ldots,5} \Vert \q_{fk} - \q_{fk}^H \Vert\,,$$
where $\q_{ik}\,, \q_{fk}$ and $\q_{ik}^H\,, \q_{fk}^H, k=0,\ldots,5$ are the B\`ezier control points of the first and second segment  of the biarc $\x$ and of the biarc representation of $\x^H,$ respectively. Since both $\x$ and $\x^H$ interpolate first order Hermite data at the end points, we have
$$\q_{i0} = \q_{i0}^H\,, \quad
\q_{i1} = \q_{i1}^H\,,\quad 
\q_{f5} = \q_{f5}^H\,,\quad 
\q_{f4} = \q_{f4}^H\,.$$ 
On the other hand, since we already know that $\x^H$ approximates $\r$ with fourth order, its second derivative at $t=0$ approximates $\r''(0)$ at least with second order. Thus, considering the expression of the second derivative in Bernstein form, this implies that it is at least $ \Vert \q_{i2} - \q_{i2}^H \Vert = O(\Delta s^4).$ 
Considering these preliminary points and also that  both biarcs $\x$ and $\x^H$ have $C^2$ smoothness at $t = 0.5,$ we can just focus on the expansions for
 \begin{equation} \label{exp}
  \Vert \q_{ik} - \q_{ik}^H \Vert, \quad k=3,4,5\,,\qquad  \Vert \q_{fk} - \q_{fk}^H \Vert, \quad k=3.
  \end{equation} 
 In order to simplify the MAPLE expansions, without loss of generality, we assume that $\frac{d\r}{ds}(0)$ and $\frac{d^2\r}{ds^2}(0)$ are aligned with the $x$ and $y$ axis, respectively, with $\frac{d^2\r}{ds^2}(0) = b_y\, {\bf j},$ and $b_y \in \R.$ Furthermore, by taking into account the arc-length parameterization for $\r,$ we set $\frac{d\r}{ds}(0) = {\bf i}$ and $\frac{d^3\r}{ds^3}(0) = -b_y^2\, {\bf i} + c_y\, {\bf j} +c_z\, {\bf k},$ with $c_y,c_z \in \R.$ The expressions of the zero and first derivative of $\r$ at $s = \Delta s$ are then symbolically defined by using Taylor expansions at $s = 0.$  
 
Now ,the expansions for the norms in (\ref{exp}) clearly depend on those of the free real parameter $a_1$, fixed using formula (\ref{a1def}), and on those of the cosine and sine of $\alpha_2$, defined through formulas (\ref{alpha2def}) and (\ref{fcoefdef}).  Considering for brevity the case $b_y \neq 0,$ 
 the symbolic implementation of the scheme produces the following formulas   
  $$ a_1 = -\frac{c_z}{16b_y} \Delta s^2 + O(\Delta s^3) \,, \; 
  \cos \alpha_2 = 1 - \frac{c_z^2}{512b_y^2} \Delta s^2 + O(\Delta s^3),\;
 \sin \alpha_2 = \frac{c_z}{16b_y} \Delta s + O(\Delta s^2)\,.$$
This implies that all the free parameters asymptotically tend to zero, in line with other PH interpolation schemes.
After the preliminary symbolic computation of the control points $\q_{ik}^H, \q_{fk}^H, k=0,\ldots,5$ defining  $\x^H$ in biarc form, the symbolic computation of all the quaternion coefficients ${\cal A}_i, {\cal B}_i, i=0,1,2$ defining the pre-image of $\x$ allows also the symbolic evaluation of $\q_{ik}, \q_{fk}, k=0,\ldots,5$ through formulas (\ref{cpi}) and (\ref{cpf}). The following expansions of the norms in (\ref{exp}) are then obtained,
\begin{align*}
&  \Vert \q_{i3} - \q_{i3}^H \Vert &= \frac{7}{(11520 \vert b_y \vert)}\sigma s^4  + O(\Delta s^5)\, ,&
& \Vert \q_{i4} - \q_{i4}^H \Vert &= \frac{1}{(1152 \vert b_y \vert)}\sigma s^4  + O(\Delta s^5)\, ,&\cr
&  \Vert \q_{i5} - \q_{i5}^H \Vert &= \frac{1}{(1152 \vert b_y \vert)}\sigma s^4  + O(\Delta s^5)\, ,&
&   \Vert \q_{f3} - \q_{f3}^H \Vert &= \frac{1}{(3840 \vert b_y \vert)}\sigma s^4  + O(\Delta s^5)\, ,&
\end{align*}
with
\begin{align*}
\sigma & = \left( b_y^8 + (24 d_y + 288 \varepsilon_y) b_y^5 + (36c_y^2 + 18c_z^2)b_y^4 - 576b_y^3c_y\varepsilon_x \right. \cr
&\qquad + (16d_y^2 + 384d_y\varepsilon_y + 16d_z^2 + 384d_z\varepsilon_z + 2304\varepsilon_x^2 + 2304\varepsilon_y^2 + 2304\varepsilon_z^2)b_y^2 \cr
&\qquad\left. -48c_z((-\frac{d_y}{2} - 6\varepsilon_y)c_z + c_y(d_z + 12\varepsilon_z))b_y + 36c_y^2c_z^2 + 9c_z^4\right) ^\frac{1}{2},
\end{align*}
where $d_y$ and $d_z$ are coefficients of the Taylor expansion of $\frac{d^4\r}{ds^4}(0)$ while $\varepsilon_x$ , $\varepsilon_y$ and $\varepsilon_z$ are the error on the second derivative, i.e. $\w_i = \frac{d^2\r}{ds^2}(0) + (\varepsilon_x , \varepsilon_y ,\varepsilon_z)^T s^2 $.

We conclude the proof reporting the asymptotic expansion of $\Vert \frac{d^2\x}{ds^2}(\Delta s) - \frac{d^2\r}{ds^2}(\Delta s) \Vert, $
 $$\Vert \frac{d^2\x}{ds^2}(\Delta s) - \frac{d^2\r}{ds^2}(\Delta s) \Vert  = \sqrt{(\varepsilon_x^2 + \varepsilon_y^2 + \varepsilon_z^2)}s^2 + O(\Delta s^3).$$
  \QED 
  
It is interesting to note that the term proportional to $s^2$ depends only on the error $\varepsilon$ on the second derivative. This means that in case of exact information, the final second derivative has a super-convergence, with $\Vert \frac{d^2\x}{ds^2}(\Delta s) - \frac{d^2\r}{ds^2}(\Delta s) \Vert_2 \,=\, O(\Delta s^3)$.
%
%

\begin{figure}[!t]
\centering
\subfigure{\hspace*{-.35cm}
\includegraphics[width=0.34\linewidth]{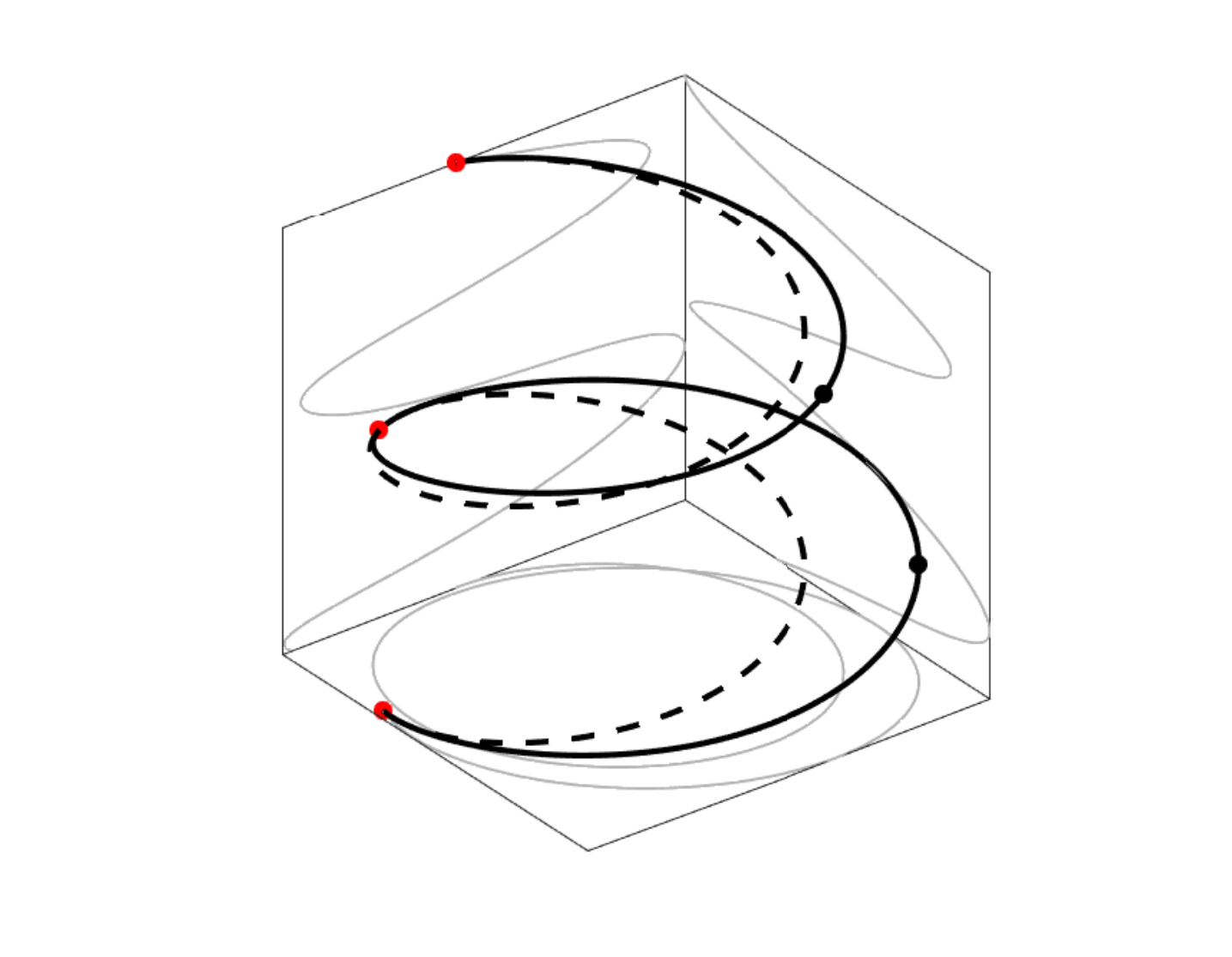}}
\subfigure{\hspace*{-.35cm}
\includegraphics[width=0.34\linewidth]{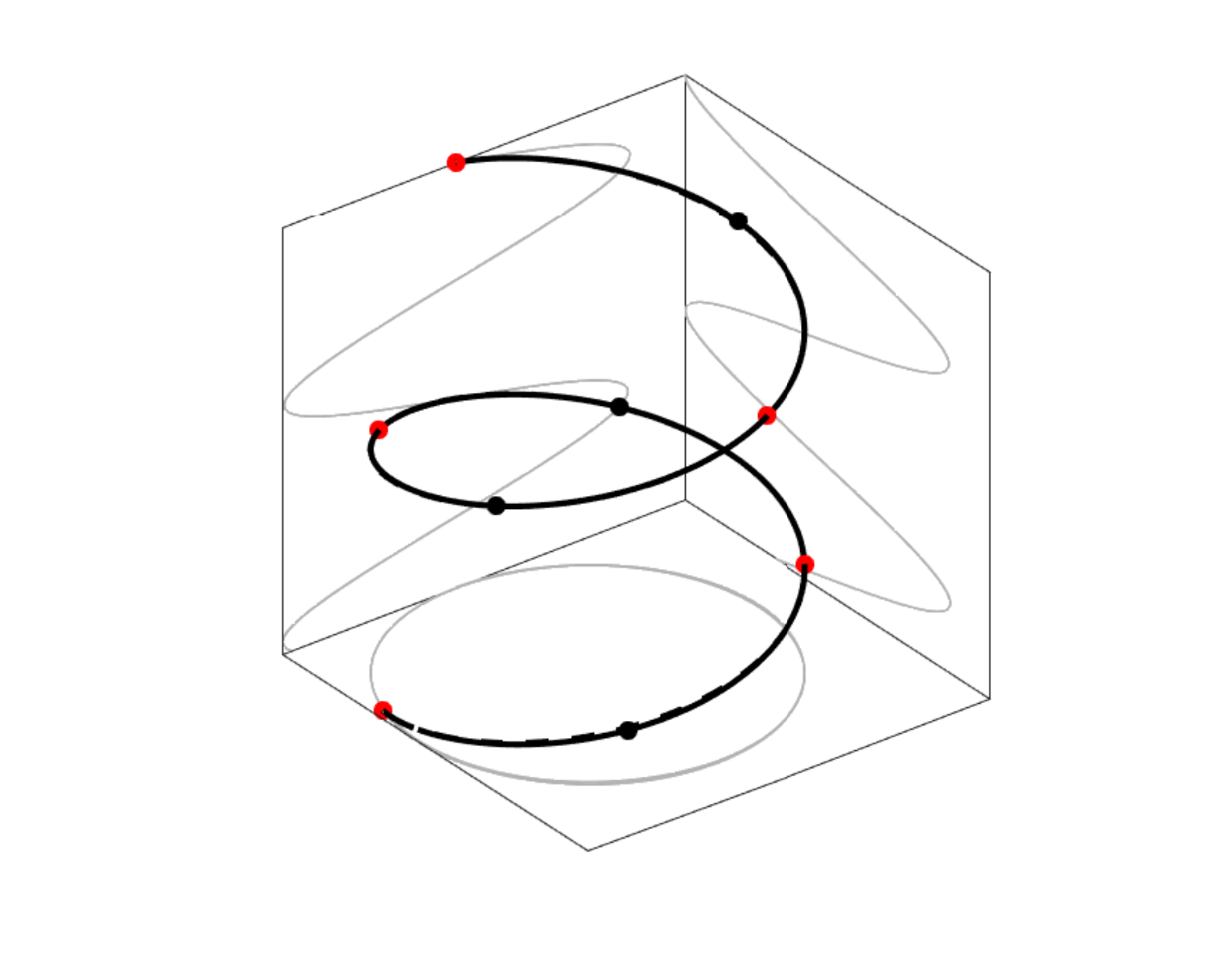}}
\subfigure{\hspace*{-.35cm}
\includegraphics[width=0.34\linewidth]{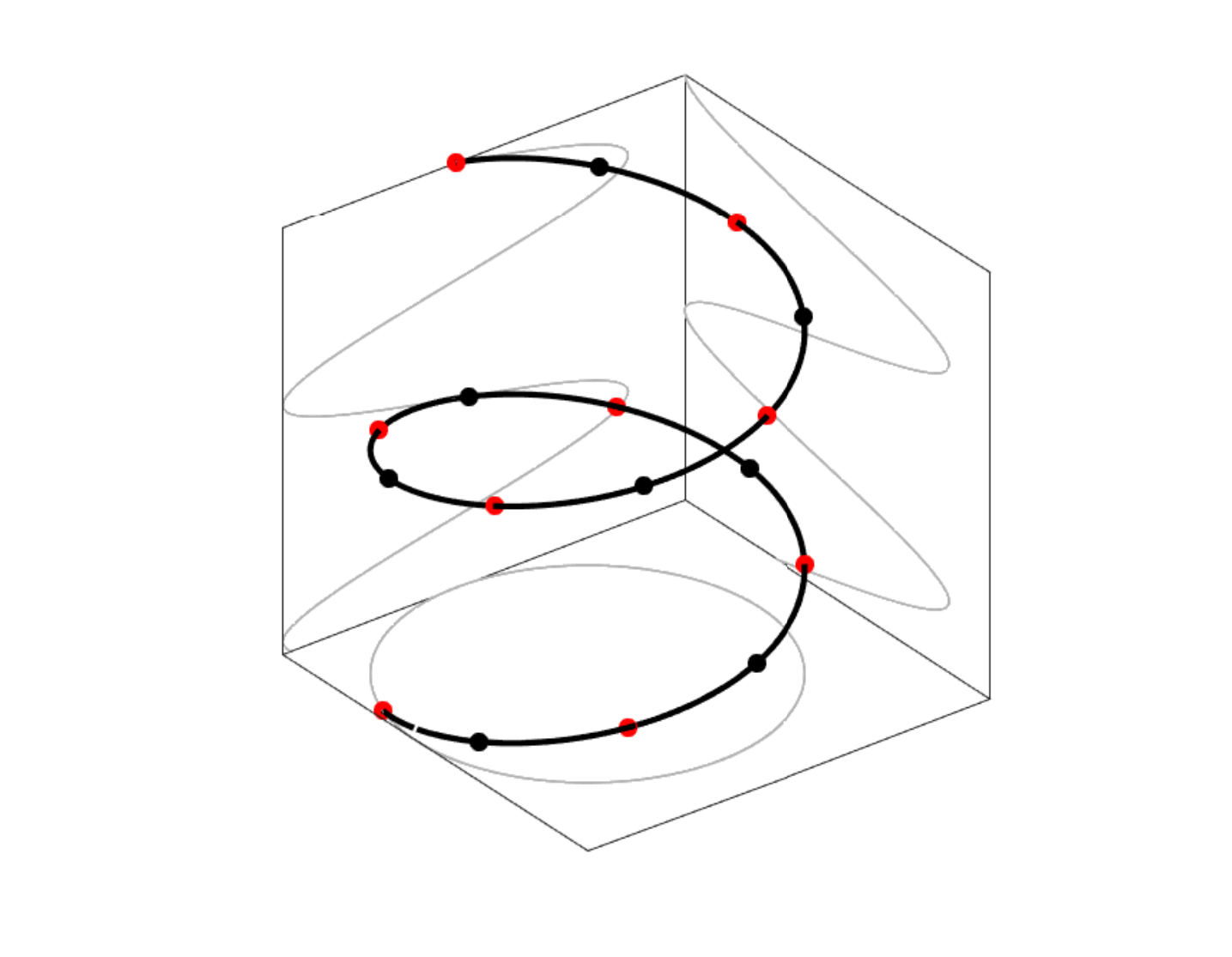}}
\caption{Reconstruction of an helix (dashed line) using the proposed algorithm, with  3 (left), 5 (center) and 9 (right) sampled points (red dots) and Hermite data. The black dots represent the joint point of two biarc segments.}
    \label{fig:helix}
\end{figure}

Figure~\ref{fig:helix} shows a first example of an helix approximation with the proposed approach. The parametric representation of the considered circular helix is given in  formula (\ref{eq:helix}) of the next section.

%% file: section05.tex
\section{Numerical results}
\label{sec:numerics}

This section presents a selection of numerical tests. Their goal is twofold. On one side, they consolidate the theoretical results obtained in the previous section concerning the approximation order of scheme. On the other side, they show its possible application to real data stream interpolation.

\subsection{Numerical approximation order}

The general procedure performed for the numerical estimation of the approximation order is now presented. We consider a generic parametric curve $\C(u)$, with $u \in [0 , U]$, to sample a set of Hermite data using the uniform parametric grid:
$$
u_j = u_j^{(k)} = j \frac{U}{N}, \, \, \mbox{ with } N = 2^k, \qquad j = 0,\ldots,N, \quad k = 0,\ldots,9.
$$
By denoting with $\X^{(k)}(u)$ the spline generated with the $C^2$ biarc algorithm so that
$$
\X^{(k)}(u) = \x_j^{(k)}(u) \qquad \mbox{for } u \in [u_j^{(k)} , u_{j+1}^{(k)}], \,\qquad j = 0,\ldots,N-1, 
$$
the interpolation problem related to each path segment $\x^{(k)}_j(u)$ is:
$$
\p_i^{(k)} = \C\left(u_j^{(k)}\right), \qquad 
\p_f^{(k)} = \C\left(u_{j+1}^{(k)}\right), \qquad 
\vv_i^{(k)} = \C'\left(u_j^{(k)}\right), \qquad 
\vv_f^{(k)} = \C'\left(u_{j+1}^{(k)}\right)
$$
with 
$$
\w_i^{(k)} = \C''\left(u_j^{(k)}\right) \quad \mbox{for } j = 0
\quad\mbox{and}\quad 
\w_i^{(k)} = \x''^{(k)}_{j-1}\left(u_j^{(k)}\right) \quad \mbox{for } j = 1,\ldots,N.
$$ 
We remind that the $'$ symbol denotes derivatives with respect to the global parameter $u$. For this test the first and second order Hermite data are obtained from certain analitic curves. We will show in the next example how to deal with point data stream. The four analitic curves considered in the numerical study are now presented.

\begin{figure}[t!]
\centering
\includegraphics[width=0.515\linewidth]{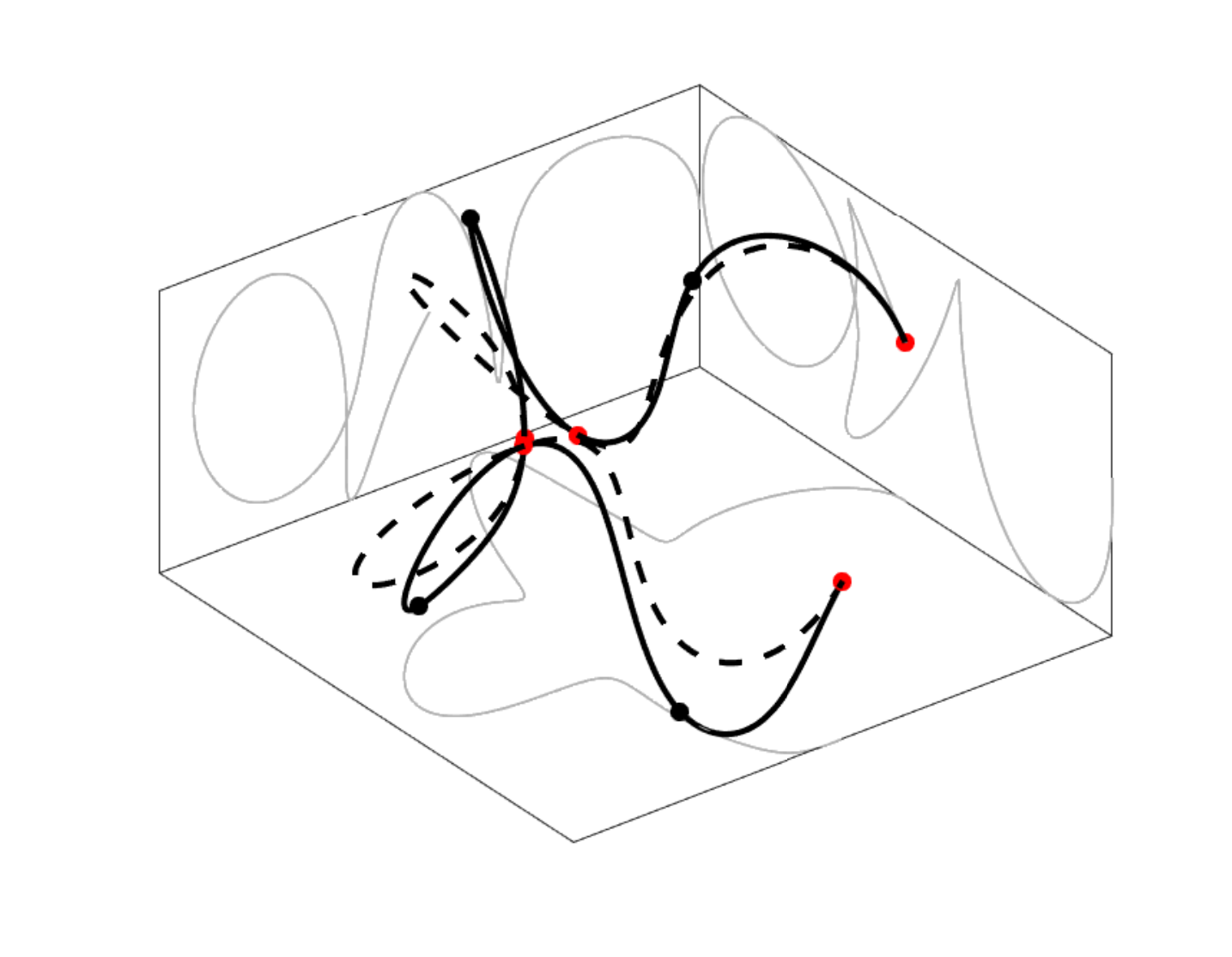}
\hspace*{-1cm}
\includegraphics[width=0.515\linewidth]{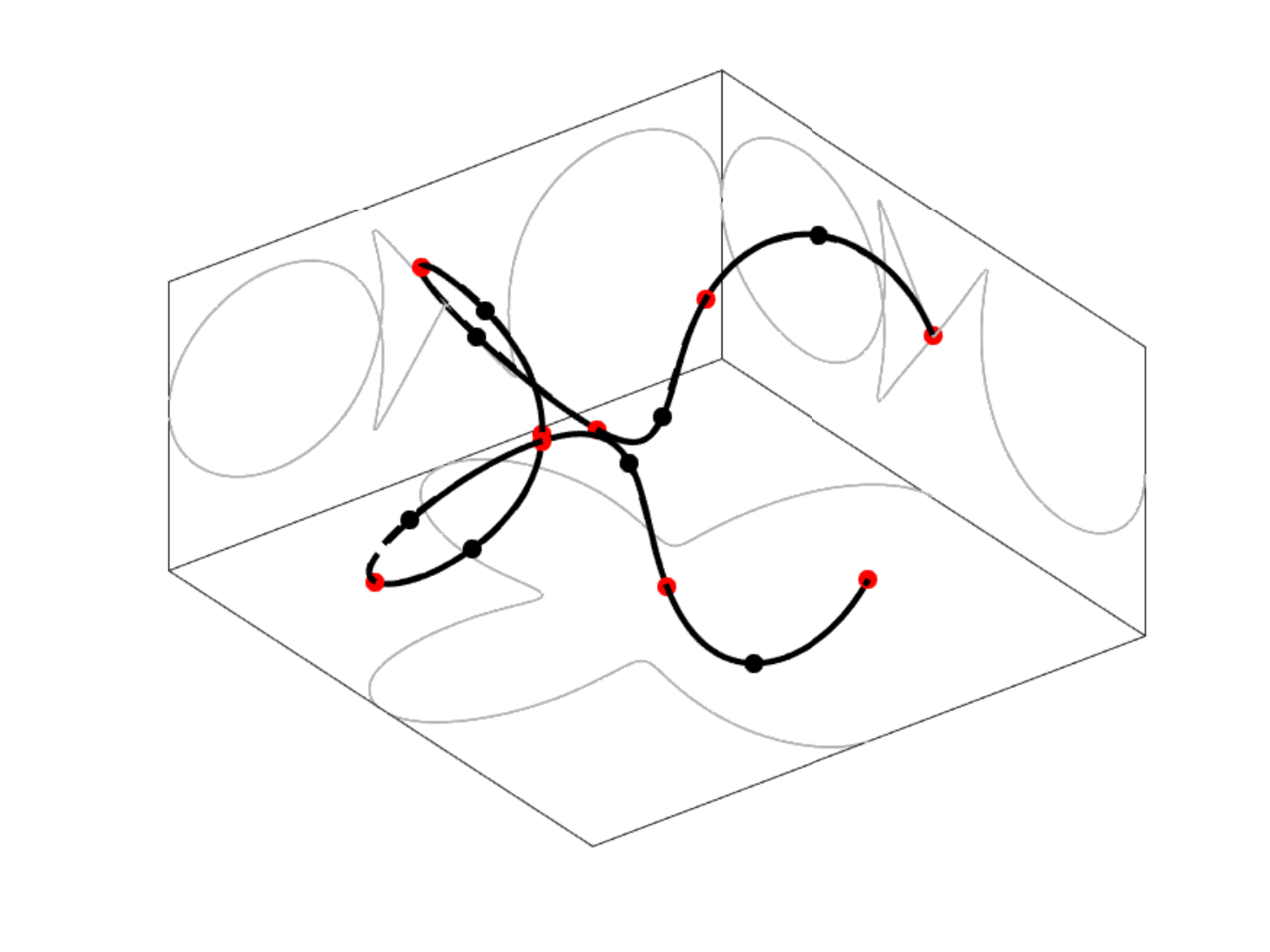}
\caption{$C^2$ PH quintic spline biarcs (solid lines) interpolating curve $\#1$ (dotted lines). The number of approximating PH quintic biarc segments is 4 on the left ($k = 2$) and 8 on the right ($k=3$). The black dots represent the joint point of two biarc segments.}
    \label{fig:curve2}
\end{figure}

\begin{figure}[t!]
\centering
\includegraphics[width=0.545\linewidth]{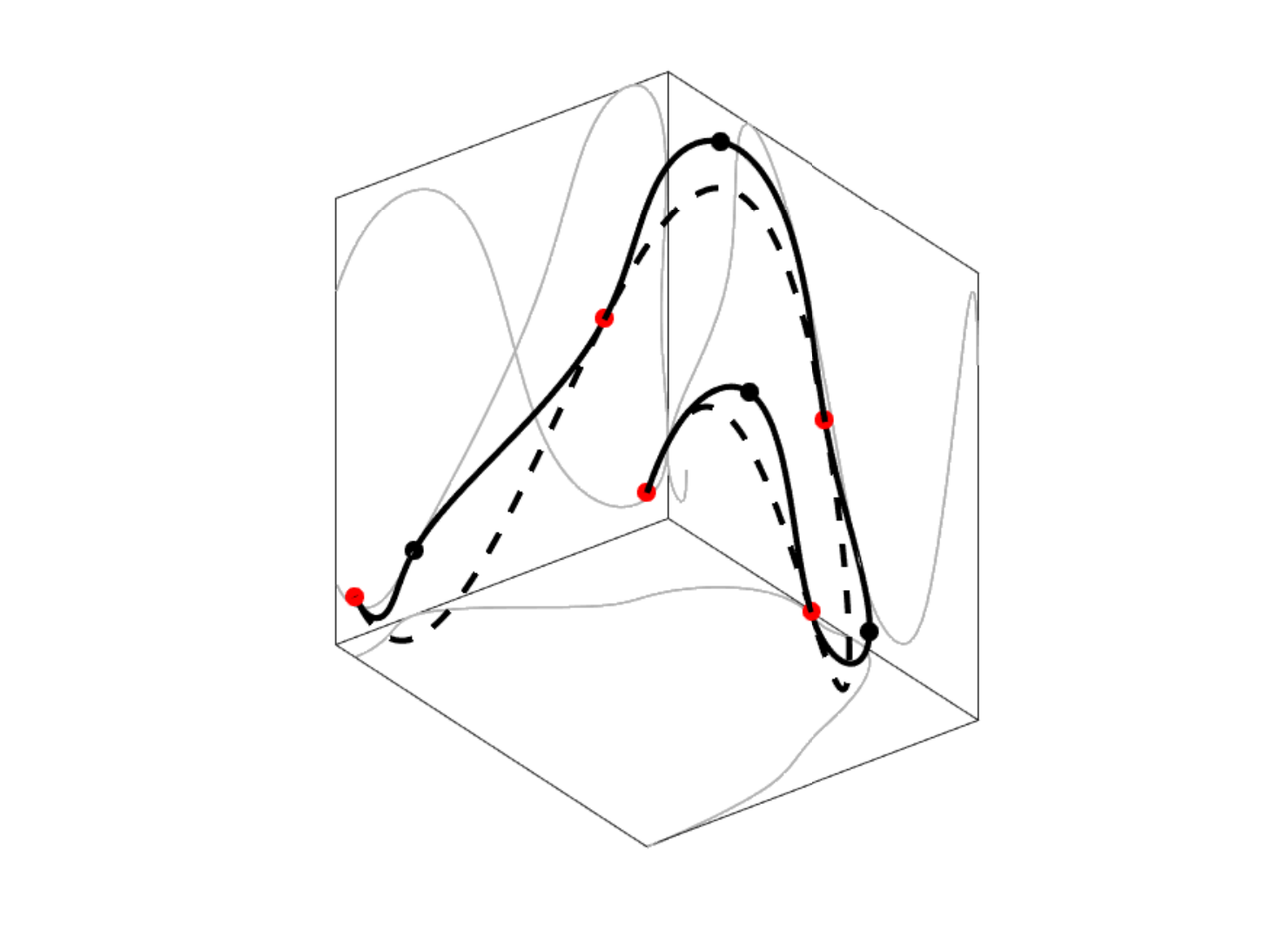}
\hspace*{-2cm}
\includegraphics[width=0.545\linewidth]{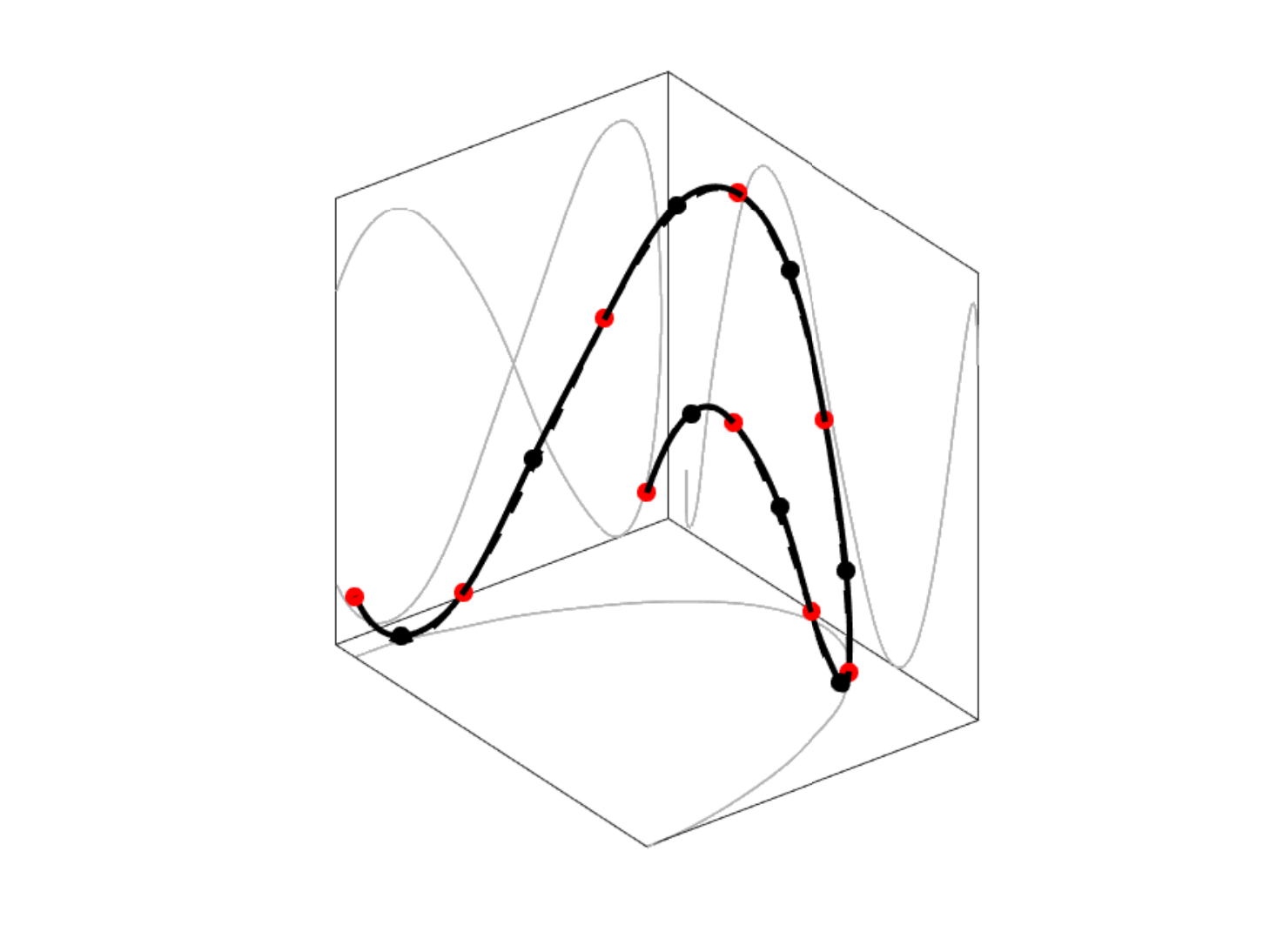}
\caption{$C^2$ PH quintic spline biarcs (solid lines) interpolating curve $\#2$ (dotted lines). The number of approximating PH quintic biarc segments is 4 on the left ($k = 2$) and 8 on the right ($k=3$). The black dots represent the joint point of two biarc segments.}
    \label{fig:curve3}
\end{figure}

\begin{figure}[t!]
\centering
\includegraphics[width=0.545\linewidth]{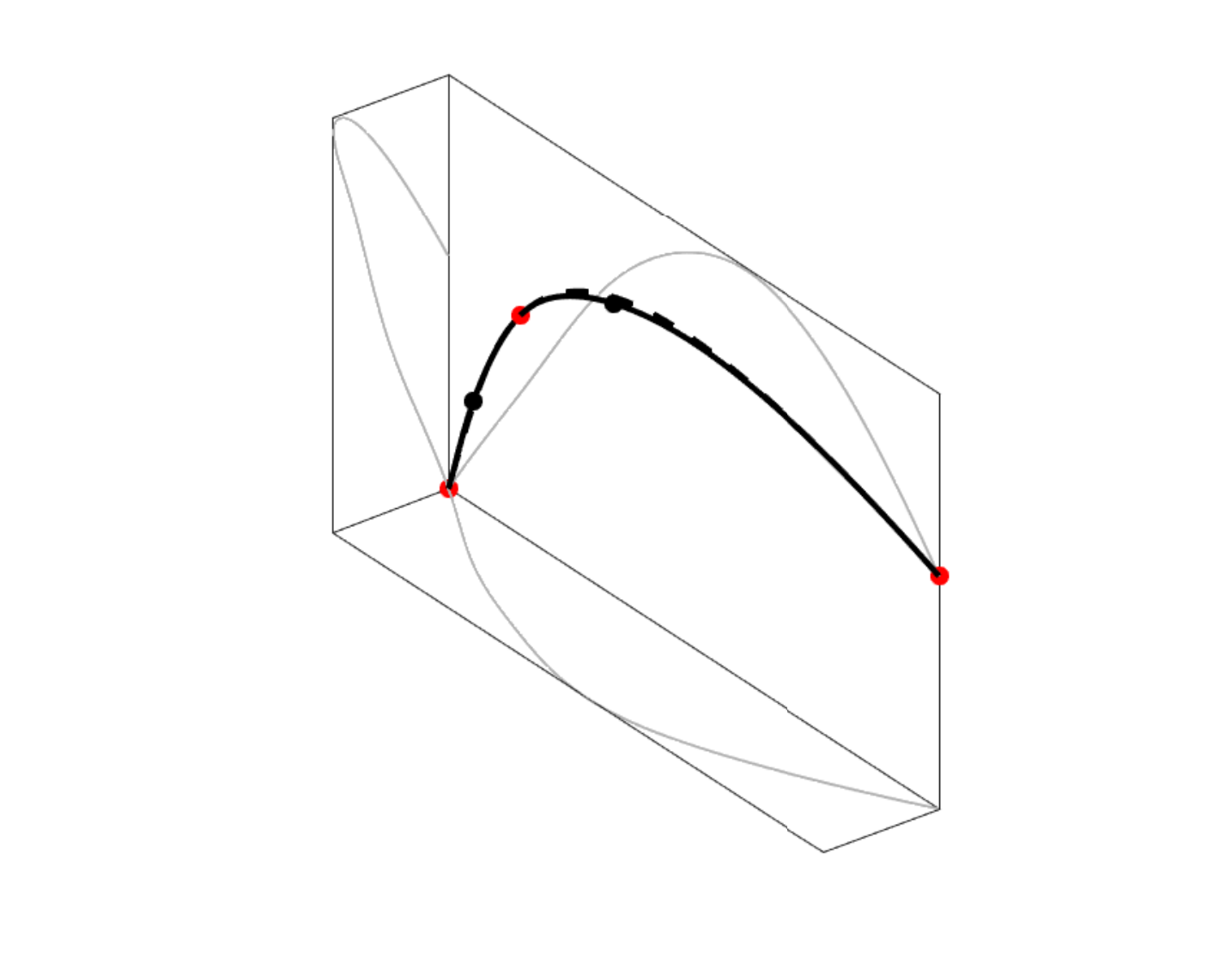}
\hspace*{-2cm}
\includegraphics[width=0.545\linewidth]{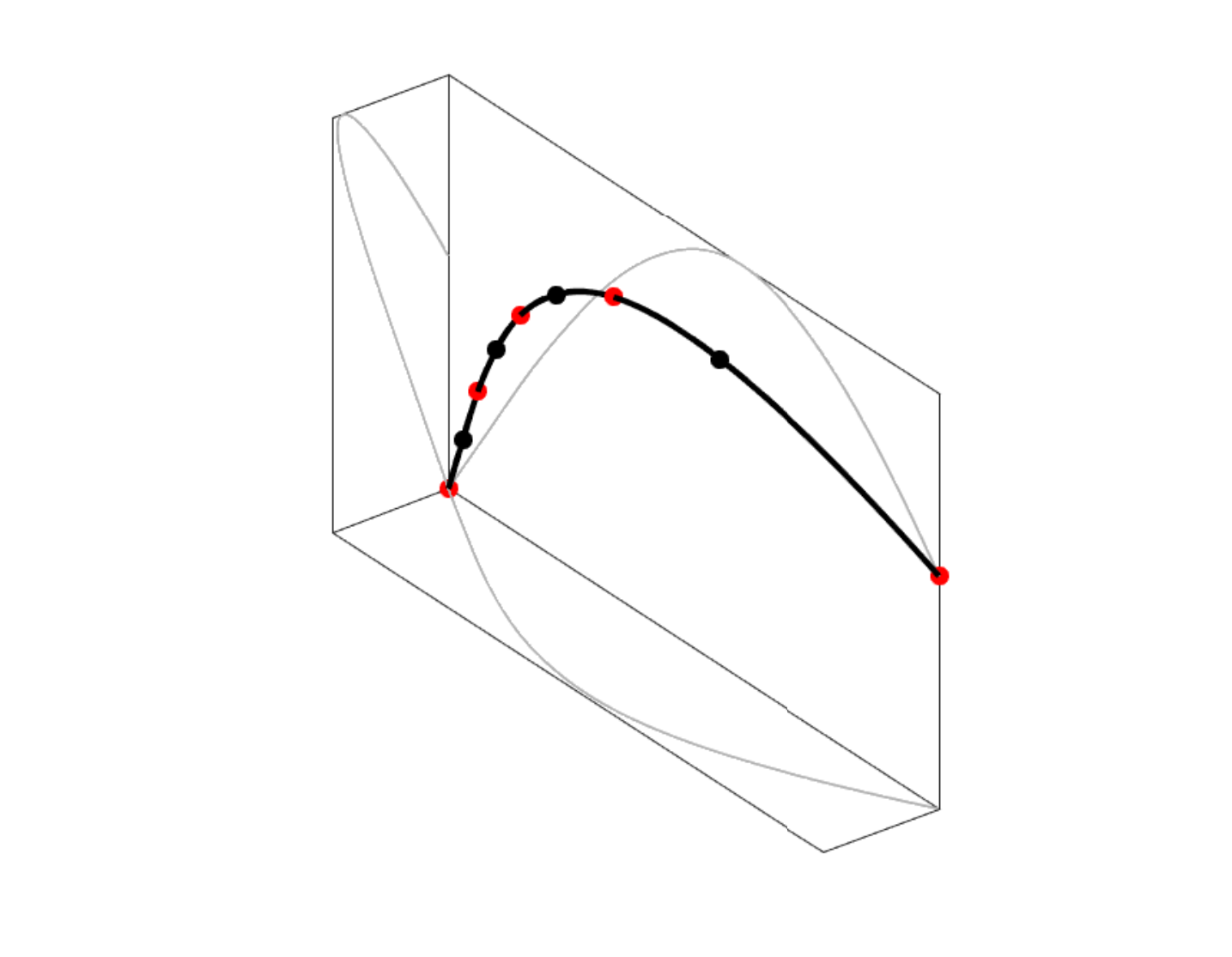}
\caption{$C^2$ PH quintic spline biarcs (solid lines) interpolating curve $\#2$ (dotted lines). The number of approximating PH quintic biarc segments is 2 on the left ($k = 1$) and 4 on the right ($k=2$). The black dots represent the joint point of two biarc segments.}
    \label{fig:curve4}
\end{figure}

In the first test (curve $\#$1) we consider data from a circular arc-length parametrized helix
\begin{equation}
\C(u) = \left(\begin{array}{l} 
10 \sin\left( \frac{u}{u_h} \right) \cr
10 \cos\left( \frac{u}{u_h} \right) \cr
-2 \frac{u}{u_h}\cr
\end{array} \right) \quad \mbox{with } u_h = \sqrt{104}, \quad u \in[0 \,,\, 3.6 \pi u_h]  
\label{eq:helix}
\end{equation}  
of total length $L = 115.34$. 
In the second test (curve $\#$2) the data are sampled from a curve on a torus with equations 
\begin{equation}
\C(u) = \left(\begin{array}{l} 
(20 + 10 \cos(3 u))\cos(0.75 u) \cr
(20 + 10 \cos(3 u))\sin(0.75 u) \cr
10 \sin(3 u)\cr
\end{array} \right), \quad u \in[0 , 2 \pi],
\label{eq:torus}
\end{equation}  
and total length $L = 212.85$. 
The Lissajous curve of the form
\begin{equation}
\C(u) = \left(\begin{array}{l} 
\cos(3 (u - u_l)) \cr
\sin(2 (u - u_l)) \cr
\sin(7 (u - u_l)) \cr
\end{array} \right), \quad \mbox{with } 
u_l = \frac{\pi}{4}, \quad u \in\left[0 , \frac{\pi}{2}\right],
\label{eq:lissajous}
\end{equation}  
and total length $L = 81.36$ is considered for the third test (curve $\#$3). Finally, we consider a curve with a zero curvature point (curve $\#$4) and equations
\begin{equation}
\C(u) = \left(\begin{array}{l} 
1 + u_z + u_z^4 - u_z^6 + u_z^8 \cr
-4 + 2u_z + u_z^5 - u_z^7 + u_z^9 \cr
2 -3u_z + u_z^3 - u_z^{10} \cr
\end{array} \right), \quad \mbox{with } u_z = \frac{u-10}{10}, \quad u \in[0 , 10].
\label{eq:zeroCurvature} 
\end{equation} 
The zero curvature point corresponds to the parameter value $u = 10$ and the total curve length is equal to $L = 4.23$. . 
Table \ref{tab:approximation} shows the approximation error 
\[
e_k = \max_{u \in [0, U]}\Vert \C(u) - X^{(k)}(u) \Vert, \quad k = 0,\ldots,9
\]
and the numerical approximation order 
$$
p_k = \log_2 \left( \frac{E_{k-1}}{E_k} \right), \quad  k = 1,\ldots,9,
$$
for the four test curves. The parameter
$u$ has been sampled at $( 2^{12} + 1 )$ uniformly spaced
 values in the interval $[0 , U]$ to compute the numerical values of $E_k$. The numerical results confirm the fourth approximation order of the interpolation scheme. Fig.~\ref{fig:helix} in the previous section shows the behavior of approximations to the circular helix (curve $\#$1) using 2, 4, and 8 PH quintic biarc interpolants. Fig.~\ref{fig:curve2}--\ref{fig:curve4} confirm that the proposed selection strategy is appropriate also for the asymptotic convergence of the other test cases (curve $\#$2, $\#$3, and $\#$4).
 

\begin{table}[!t]
\begin{tabular}{ p{0.3cm} | p{2cm} p{1.1cm}| p{2cm} p{1.1cm} | p{2cm} p{1.1cm}| p{2cm} p{1.1cm} }
 \hline
  & \multicolumn{2}{|c|}{curve $\#$1} & \multicolumn{2}{|c|}{curve $\#$2} &\multicolumn{2}{|c|}{curve $\#$3} &\multicolumn{2}{|c}{curve $\#$4}\\ 
 \hline
 $k$ & $e_k$ & $p_k$ &$e_k$ & $p_k$ &$e_k$ & $p_k$ &$e_k$ & $p_k$\\
 \hline
 \hline
 2 & 1.6891e-01 & 5.00 & 1.1459e+01 & 3.76 & 2.7103e+00 & 2.78 & 3.1454e-03 & 3.61 \\
 3 & 1.2229e-02 & 3.79 & 2.1602e-01 & 2.77 & 3.9744e-01 & 5.73 & 1.6078e-04 & 4.29 \\
 4 & 1.0656e-03 & 3.52 & 2.3278e-02 & 3.32 & 3.9717e-02 & 3.21 & 1.1994e-05 & 3.74 \\
 5 & 8.0828e-05 & 3.72 & 2.6217e-03 & 4.42 & 1.8504e-03 & 3.15 & 8.3209e-07 & 3.85 \\
 6 & 5.5547e-06 & 3.86 & 2.1238e-04 & 3.01 & 2.3004e-04 & 3.63 & 5.5015e-08 & 3.92 \\
 7 & 3.6361e-07 & 3.93 & 1.4952e-05 & 3.78 & 1.6709e-05 & 3.83 & 3.5399e-09 & 3.96 \\
 8 & 2.3249e-08 & 3.97 & 9.8900e-07 & 3.94 & 1.0883e-06 & 3.92 & 2.9747e-10 & 3.57 \\
 9 & 1.4695e-09 & 3.98 & 6.3543e-08 & 3.98 & 6.8957e-08 & 3.96 & 1.8632e-11 & 4.00 \\
 \hline
\end{tabular}
\caption{Approximation error $e_k$ and related numerical approximation order $p_k$ for the four curves described by equations \eqref{eq:helix}, \eqref{eq:torus}, \eqref{eq:lissajous} and \eqref{eq:zeroCurvature}.}
\label{tab:approximation}
\end{table}


\subsection{Application to 3D point stream interpolation}

We now present a test for the application of the $C^2$ PH biarc interpolation algorithm to 3D data stream interpolation by considering a sequence of points $\p_j$ for $j = 0,\ldots,N$. As first step, it is necessary to set a specific global parameterization. The test is based on the chord-length parametization, that allows to choose the parameter value associated to each point as
$$
u_j = u_{j-1} - \Vert \p_j - \p_{j-1}  \Vert
$$ 
with $u_0= 0$. By considering only input point streams, local rules for derivative approximations has to be coupled with the inteprolation scheme.
We rely on the local formulas called MinAJ2 introduced in \cite{Debski1} for data stream application. They enable the construction of suitable PH spline interpolants with fair shapes, as already mentioned by the author in the context of standard splines. These formulas have also a simple implementation since during the stream elaboration the tangent direction is available at the left point, while the right one is simply defined as the first derivative of a standard local $C^2$ cubic spline. When constructing an inner spline section, the local cubic spline is obtained by requiring the interpolation of the two biarc end-points, together with the left tangent and the successive stream point, suitably combined with the minimization of a suitable fairness functional. The Hermite data are then computed with a short delay, equivalent to the time necessary for one point stream acquisition. More precisely, the right derivative $\vv_j$ is chosen as
$$
\vv_j = \frac{A\,\p_{j-1} + B\,\vv_{j-1} + C\,\p_j + D\,\p_{j+1}}{E}, \qquad j = 1,\ldots,N-1,
$$
where
$$
\left\{ \begin{array}{l} 
A = - (u_{j+1} - u_j)^2(2 u^2_{j+1} + 2u_ju_{j+1} - u^2_j), \cr
B = - u_ju^2_{j+1}(u_{j+1} - u_j)^2, \cr
C = u_{j+1}(2 u^3_{j+1} - 2u_ju^2_{j+1} - 3u^2_ju_{j+1} + u^3_j), \cr
D = u^3_j(2u_{j+1} - u_j), \cr
E = - u_ju_{j+1}(u_{j+1} - u_j)(2 u^2_{j+1} + 2u_ju_{j+1} - u^2_j), \cr
\end{array} \right.
$$
The first and the last derivative are computed as 
$$
\vv_0 = \frac{(\p_1 - \p_0)(u_2 - u_0)^2 + (\p_1 - \p_2)(u_1 - u_0)^2}{(u_1 - u_0)(u_2 - u_0)(u_2 - u_1)}
$$
and
$$
\vv_{N} = -\frac{(\vv_{N-1}(u_N-u_{N-1}) - 2\p_N + 2\p_{N-1})}{(u_N-u_{N-1})},
$$
respectively.
We refer to \cite{Debski1} for further details. Finally, in order to avoid the choice of an arbitrary value for $\w_0$, the first segment of the spline is computed with the CC PH algorithm presented in \cite{farouki08c}.

\begin{figure}[!t]
\centering
\includegraphics[width=0.45\linewidth]{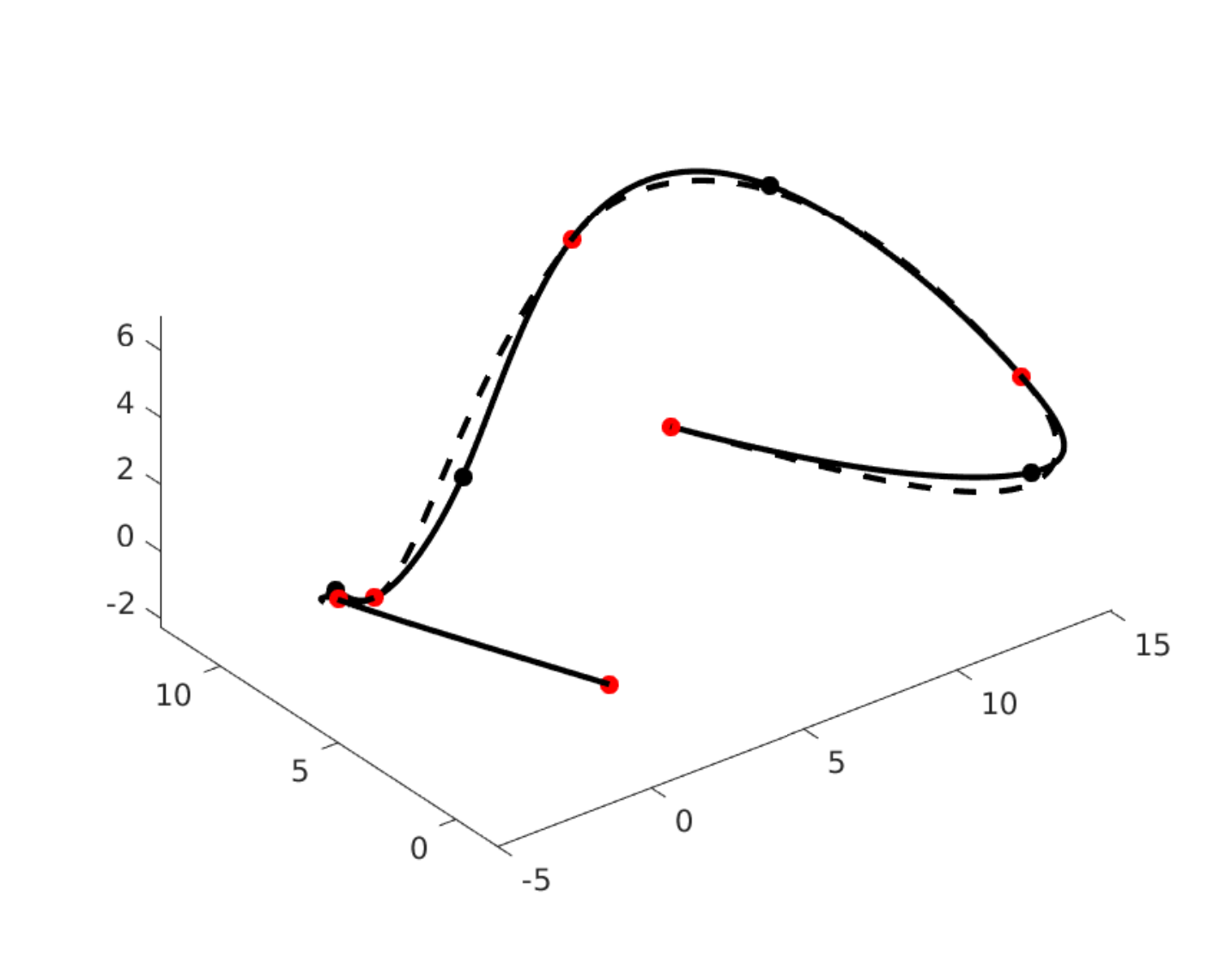}
\includegraphics[width=0.45\linewidth]{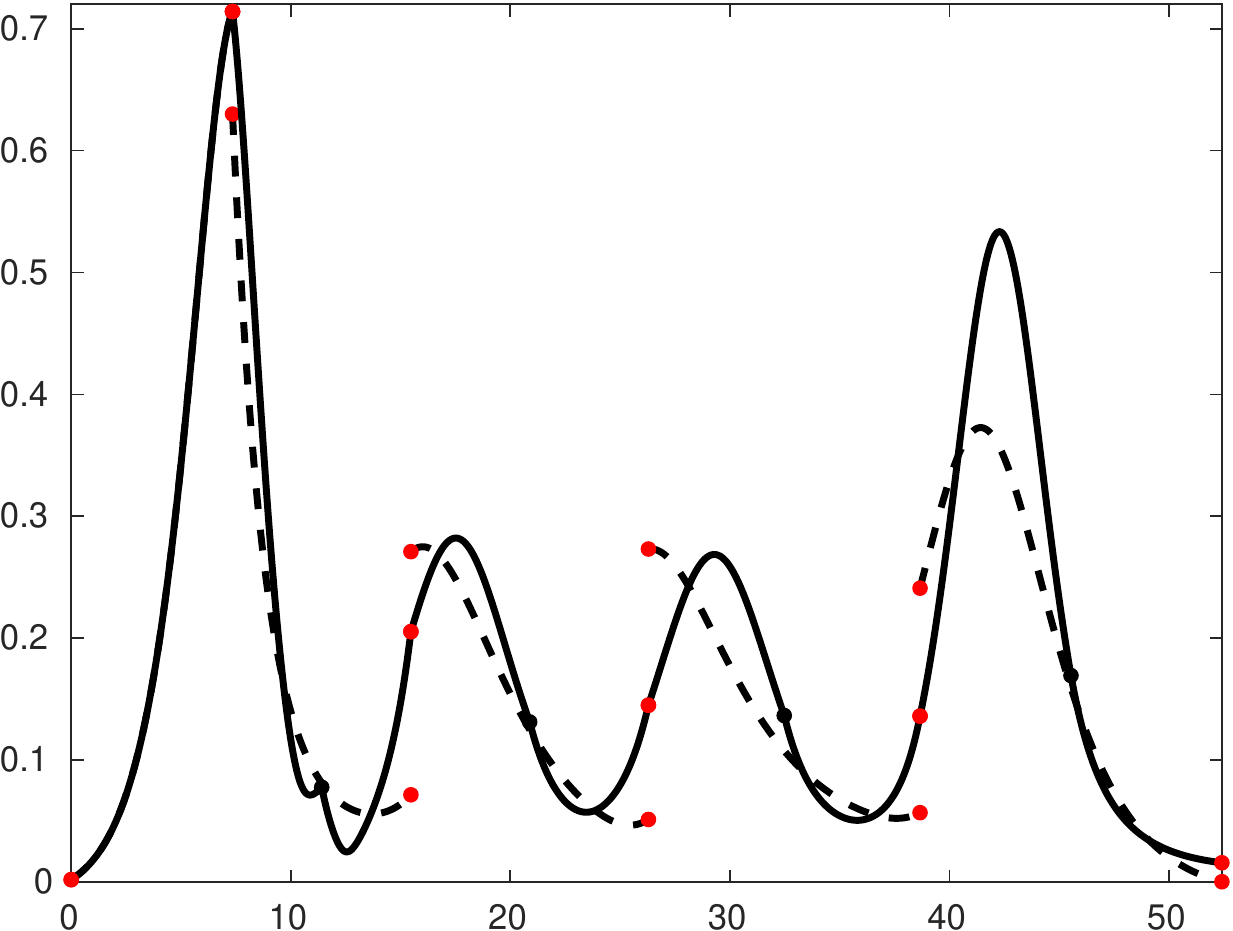}
\caption{
$C^2$ PH quintic spline biarcs (solid lines, left) and $C^1$ CC PH quintic spline (dashed lines, left) interpolating a point sequqnce (red dots, left). The black dots represent the joint point of two biarc segments. The curvature plots of the two spline paths are also shown (right).}
    \label{fig:DataStream}
\end{figure}

The 3D data stream used in the test is the following:
$$
\begin{array}{lll} 
\p_0 = (0 , 0 , 0)^T, & \p_0 = (-5 , 5 , 2)^T, & \p_0 = (0 , 10 , -2)^T, \cr
\p_0 = (8 , 12 , 5)^T, & \p_0 = (15 , 2 , 3)^T, & \p_0 = (2 , 0 , 7)^T. \cr
\end{array}
$$ 
Figure~\ref{fig:DataStream} shows the $C^2$ PH quintic spline interpolant obtained with the biarc construction here proposed and its curvature plot, together with the comparison with the $C^1$ PH quintic spline obtained by considering the CC selection strategy. It is clear that our new PH spline construction preserves very nice smoothness property, while simultaneously ensuring the appealing feature of $C^2$ continuity.  

%% file: appendix.tex
\section*{Appendix}\label{sec:appendix}

In this appendix we recall   the basic rules of the non commutative quaternion algebra $\mathbb{H},$ used in the paper. Each quaternion ${\cal Q} \in \mathbb{H}$ can be defined as $(q_0,q_1,q_2,q_3)^T\,,$ with $q_i \in \R,$ and with $q_0$ and $\q :=(q_1,q_2,q_3)^T$ respectively referred to as {\it scalar} and {\it vector}  part of the quaternion ${\cal Q}.$ With this notation a short scalar/vector representation can also be adopted for ${\cal Q},$

$${\cal Q} = q_0 + \q\,,$$ 
where, if $q_0 = 0\,, {\cal Q}$ is said a {\it pure vector} quaternion and can be shortly denoted just as $\q.$ Conversely, when $\q$ vanishes, ${\cal Q}$ is a {\it pure scalar} quaternion and can just be denoted as any real number.  
The quaternion sum in $\mathbb{H}$ is the standard sum in $\R^4$ but the quaternion product has a specific non commutative definition that can be compactly defined as 

$$ {\cal A} {\cal B} = (a_0 + \a)(b_0+\b) = (a_0b_0 - \a \cdot  \b) + (a_0 \b + b_0 \a + \a \times \b)\,,$$
where standard notation to denote scalar and cross vector products is used. 
The {\it conjugate} of a quaternion ${\cal Q}$ is denoted as ${\cal Q}^*$ and defined as ${\cal Q}^* := q_0 - \q.$ This implies   that ${\cal Q}{\cal Q}^* = {\cal Q}^*{\cal Q} =q_0^2 + \q^T\q$ is just a pure scalar quaternion. We also observe  that, for any vector $\vv$ and quaternion  ${\cal Q}$, the quaternion product of the form ${\cal Q}\,\vv\, {\cal Q}^*$  defines a pure vector quaternion. The {\it module} $\vert {\cal Q} \vert$ of a quaternion is defined as $\vert {\cal Q} \vert := \sqrt{{\cal Q}{\cal Q}^*}$ and ${\cal Q}$ is a {\it unit} quaternion if $\vert {\cal Q} \vert = 1.$ Unit quaternions allow a compact representation of spatial rotations. For any pure vector quaternion $\vv$ and unit quaternion
${\cal Q} = \cos({\theta}/{2}) + \w \sin({\theta}/{2})$, the product ${\cal Q}\, \vv\, {\cal Q}^*\,,$ always defines a pure vector quaternion, that corresponds to a rotation of $\vv$ through angle $\theta$ about the axis defined by $\w$.